\pdfoutput=1
\documentclass[final]{siamltex}

\usepackage[scaled=.92]{helvet}
\usepackage{times}
\usepackage{booktabs}
\usepackage{wrapfig}
\usepackage{graphicx}
\usepackage{eucal} 
\usepackage{amsmath,amssymb,bm}
\usepackage{array}
\usepackage{multirow}
\usepackage{subfigure}
\usepackage{tikz}
\usepackage{pgfplots}
\usepackage{cprotect}
\usepackage{comment}
\usepackage{hyperref}
\usepackage{arydshln}

\usepackage[ruled,vlined,linesnumbered]{algorithm2e}

\newcommand{\bx}{\mathbf{x}}    
\newcommand{\bu}{\mathbf{u}}    
\newcommand{\bb}{\mathbf{b}}    
\newcommand{\bv}{\mathbf{v}}

\newcommand{\bg}{\mathbf{g}}
\newcommand{\bbf}{\mathbf{f}}
\newcommand{\bK}{\mathbf{K}}
\newcommand{\bM}{\mathbf{M}}
\newcommand{\bA}{\mathbf{A}}
\newcommand{\bB}{\mathbf{B}}
\newcommand{\bC}{\mathbf{C}}
\newcommand{\tbA}{\widetilde{\mathbf{A}}}
\newcommand{\bL}{\mathbf{L}}

\newcommand{\bI}{\mathbf{I}}

\newcommand{\bU}{\mathbf{U}}
\newcommand{\bV}{\mathbf{V}}
\newcommand{\bX}{\mathbf{X}}
\newcommand{\bG}{\mathbf{G}}
\newcommand{\bW}{\mathbf{W}}
\newcommand{\bT}{\mathbf{T}}

\newcommand{\bP}{\mathbf{P}}
\newcommand{\bS}{\mathbf{S}}

\newcommand{\bzero}{\mathbf{0}}
\newcommand{\sU}{\bm{\mathcal{U}}}
\newcommand{\sC}{\mathcal{C}}
\newcommand{\sR}{\mathcal{R}}

\newcommand{\bbR}{\mathbb{R}}
\newcommand{\bbD}{\mathbb{D}}
\newcommand{\bbB}{\mathbb{B}}
\newcommand{\bbZ}{\mathbb{Z}}
\newcommand{\bbL}{\mathbb{L}}

\newcommand{\rDkToj}{R^{D}_{k \rightarrow j}}
\newcommand{\rDkTojs}{R^{D}_{k \rightarrow j,s}}
\newcommand{\rOkToj}{R^{O}_{k \rightarrow j}}
\newcommand{\srDkToj}{\sR^{D}_{k \rightarrow j}}
\newcommand{\srDkTojs}{\sR^{D}_{k \rightarrow j,s}}
\newcommand{\srOkToj}{\sR^{O}_{k \rightarrow j}}

\newcommand{\bLoj}{\bL^{O}_{j}}
\newcommand{\sUoj}{\sU^{O}_{j}}
\newcommand{\bLok}{\bL^{O}_{k}}
\newcommand{\bLdj}{\bL^{D}_{j}}
\newcommand{\dBlock}[1]{\bL^{D}_{j, #1}}
\newcommand{\schurBlock}[1]{\sU^{D}_{j, #1}}
\newcommand{\uBlock}[1]{\bU^{D}_{j, #1}}
\newcommand{\vBlock}[1]{\bV^{D}_{j, #1}}
\newcommand{\bAdj}{\bA^{D}_{j}}
\newcommand{\sUdj}{\sU^{D}_{j}}
\newcommand{\bLdk}{\bL^{D}_{k}}
\newcommand{\bLbi}{\bL^{\bbB}_{i}}

\newcommand{\rDs}{R^{D}_{j,s}}
\newcommand{\srDs}{\sR^{D}_{j,s}}
\newcommand{\rD}[1]{R^{D}_{j, #1}}
\newcommand{\cDs}{C^{D}_{j,s}}
\newcommand{\scDs}{\sC^{D}_{j,s}}
\newcommand{\cD}[1]{C^{D}_{j, #1}}

\newcommand{\sCbi}{\sC^{\bbB}_{i}}
\newcommand{\sRbi}{\sR^{\bbB}_{i}}

\newcommand{\bUj}{\bU_{j}}
\newcommand{\bVj}{\bV_{j}}
\newcommand{\bUk}{\bU_{k}}
\newcommand{\bVk}{\bV_{k}}

\DeclareMathOperator*{\blockRank}{rank}
\DeclareMathOperator*{\chol}{chol}
\DeclareMathOperator*{\scatter}{scatter}
\DeclareMathOperator*{\gatherRows}{gatherRows}
\DeclareMathOperator*{\scatterRows}{scatterRows}
\DeclareMathOperator*{\scatterColumns}{scatterColumns}
\DeclareMathOperator*{\alignSet}{alignSet}
\DeclareMathOperator*{\cholesky}{cholesky}
\DeclareMathOperator*{\factorSupernode}{factorSupernode}
\DeclareMathOperator*{\offDiagonalMultiply}{offDiagonalMultiply}
\DeclareMathOperator*{\offDiagonalMultiplyTrans}{offDiagonalMultiplyTranspose}
\DeclareMathOperator*{\diagonalMultiply}{diagonalMultiply}
\DeclareMathOperator*{\diagonalMultiplyTrans}{diagonalMultiplyTranspose}
\DeclareMathOperator*{\approximateOffDiagonal}{approximateOffDiagonal}
\DeclareMathOperator*{\approximateDiagonalBlock}{approximateDiagonalBlock}
\DeclareMathOperator*{\factorDiagonal}{factorDiagonal}
\DeclareMathOperator*{\diagonalSolve}{diagonalSolve}
\DeclareMathOperator*{\randomMatrix}{randomMatrix}
\DeclareMathOperator*{\makeOrthonormal}{makeOrthonormal}

\newcommand{\figspace}{\vspace{-0.4cm}}

\newcommand{\myparagraph}[1]{\noindent \textsf{\textbf{#1}}}

\newlength\savedwidth
\newcommand\whline[1]{\noalign{\global\savedwidth\arrayrulewidth 
                               \global\arrayrulewidth #1} %
                      \hline
                      \noalign{\global\arrayrulewidth\savedwidth}}

\newcommand{\todo}[1]{{$\spadesuit$ {\bf #1} $\spadesuit$ }}

\title{An Efficient Solver for Sparse
Linear Systems Based on Rank-Structured Cholesky Factorization} 

\author{Jeffrey N. Chadwick \and David S. Bindel}

\begin{document}
\maketitle

\begin{abstract}
  Direct factorization methods for the solution of large, sparse
  linear systems that arise from PDE discretizations are robust, but
  typically show poor time and memory scalability for large systems.
  In this paper, we describe an efficient sparse, rank-structured
  Cholesky algorithm for solution of the positive definite linear
  system $\bA \bx = \bb$ when $\bA$ comes from a discretized
  partial-differential equation.  Our approach combines the efficient
  memory access patterns of conventional supernodal Cholesky
  algorithms with the memory efficiency of rank-structured direct
  solvers.  For several test problems arising from PDE
  discretizations, our method takes less memory than standard sparse
  Cholesky solvers and less wall-clock time than standard
  preconditioned iterations.
\end{abstract}



\pagestyle{myheadings}
\thispagestyle{plain}
\markboth{JEFFREY N. CHADWICK AND DAVID S. BINDEL}
{AN EFFICIENT SOLVER BASED ON RANK-STRUCTURED
CHOLESKY FACTORIZATION}

%
%
%
%

\section{Introduction}

We consider the problem of solving a sparse linear system
\begin{equation}
\label{eq:main}
\bA \bx = \bb
\end{equation}
in which $\bA$ is symmetric and positive definite (SPD).
In particular, we consider the Cholesky factorization 
\begin{equation}
\label{eq:cholesky}
  \bA = \bL \bL^T,
\end{equation}
where $\bL$ is a sparse lower triangular matrix; 
see~\cite{george81,duff86,davis06}.
After factoring $\bA$, one solves~\eqref{eq:main} by two triangular
solves, at cost proportional to the number of nonzeros in $\bL$.
This approach solves~\eqref{eq:main} exactly up to roundoff effects,
and modern supernodal factorization algorithms achieve high flop rates
by organizing the factorization around dense matrix kernels.
The chief drawback of sparse direct methods is that the factor $\bL$
may generally have many more nonzero elements than $\bA$.
These {\em fill} elements limit scalability of the method in both
time and memory used, particularly for problems coming from
the discretization of three-dimensional PDEs, where the number of
nonzeros in $\bL$ typically scales as $O(N^{3/2})$, where $N$ is the
dimension of $\bA$.

Compared to direct factorization, iterative methods
for~\eqref{eq:main} generally cost less in memory and in time per step than
direct methods, but converge slowly without a good preconditioner.
Preconditioning involves a complex balance between the progress in
each step and the cost of setting up and applying the preconditioner.
Even for a single preconditioner type, there are usually many
parameters that are optimized on a problem-by-problem basis.
For this reason, packages like PETS provide interfaces to
allow users to quickly experiment with different preconditioners and
parameter settings~\cite{petsc-user-ref}, while commercial finite element codes
often forego the potential benefits of iterative methods and simply
use out-of-core direct solvers~\cite{poole2003advancing}.

Fast direct factorization methods and preconditioned iterative solvers
each use a different types of structure.
A key idea behind sparse direct methods is that one can use the
structure of the graph associated with $\bA$ to reason about
about fill in $\bL$.
This graph-theoretic approach underlies many modern sparse matrix
algorithms, from methods of computing fill-reducing elimination
orderings to ``supernodal'' factorization methods organized
around dense matrix operations on columns with similar nonzero
structure~\cite{davis06}.
In contrast, to solve problems arising from elliptic PDE
discretization efficiently, multi-level preconditioners exploit
the elliptic regularity of the underlying differential equation.
Building on ideas from fast direct solvers for integral
equations~\cite{greengard09}, recent work in ``data-sparse'' direct
solvers uses both types of structure at once, computing a structured
factorization that incorporates (approximate) low-rank
blocks~\cite{grasedyck2008parallel,xia2009,martinsson09,schmitz2014fast,gillman2014direct}.

In this paper, we describe an efficient sparse, rank-structured
Cholesky algorithm for the solution of~\eqref{eq:main} when $\bA$
comes from discretization of a PDE.  Our method combines the efficient
memory access patterns of conventional supernodal Cholesky algorithms
with rank-structured direct solvers.  Unlike prior solvers, our method
works as a ``black box'' solver, and does not require information
about
an underlying PDE mesh.
For several test problems
arising from PDE discretizations, we show that our method takes less
memory than standard sparse Cholesky codes and wall-clock time than
standard preconditioners
The remainder of the paper is organized as follows.  
In Section~\ref{SECbackground}, we briefly review the standard supernodal
left-looking sparse Cholesky algorithm on which our method is based.
In the supernodal factorization, each supernode has an associated
diagonal block storing interactions within that supernode, and an
off-diagonal block for storing interactions between supernodes.
In Section~\ref{SECoffdiagcomp}, we describe how our algorithm forms
and uses low-rank approximations to the off-diagonal blocks of the
supernodes, 
and in Section~\ref{SECdiagcomp}, we describe our approach to
hierarchical compression of the diagonal blocks.
We discuss some key implementation details in
Section~\ref{SECoptimizations},
and illustrate the behavior of our algorithm on several
example problems in~\ref{SECresults}. 
Finally, in Section~\ref{SECconclusions}, we conclude
and give potential directions for future work.

\section{Background and Notation}
\label{SECbackground}

We focus primarily on {\em supernodal left-looking Cholesky
factorization}.  This method has yielded implementations which make effective
use of modern computing architectures to efficiently solve \eqref{eq:main}
\cite{chen08}.

\subsection{Supernodal Left-Looking Cholesky Factorization}
\label{SECsupernodal}

Most sparse Cholesky codes have two phases: a fast symbolic analysis
phase to compute the nonzero structure of $\bL$, and a more expensive
numerical factorization phase in which the actual elements of $\bL$
are computed.
The symbolic analysis phase is organized around an {\em elimination
  tree} that encodes the structure of $\bL$: in general, 
$l_{ij} \neq 0$ precisely when there is some $k$ such that 
$a_{ik} \neq 0$ and $j$ is reachable from $k$ by an elimination 
tree path that passes only through nodes with indices less than $i$.
Often, the elimination tree has chains of sequentially-numbered nodes
corresponding to columns with similar nonzero structure;
these can be seen as {\em supernodes} in a coarsened version of the
elimination tree.
Supernodal factorization algorithms organize $\bL$ around such
supernodes, formed by collecting adjacent columns which are
predicted to have similar non-zero patterns in $\bL$. 
Supernodal methods achieve high efficiency by storing the nonzero
entries for a supernode together as a dense matrix, and by operating
on that matrix with optimized kernels from the BLAS.

Suppose $\bL \in \bbR^{N \times N}$ is partitioned into $M$ supernodes.  
Let $\sC_j = \left( c : c_{j} \leq c < c_{j+1} \right)$ refer to the
column indices in supernode $j$, and let 
$\sC^O_j = \left( c : c \geq c_{j+1} \right)$ be the list of columns 
occurring {\em after} supernode $j$.
Finally, let $\bL_{j}$ refer to the
block column $\bL_{j} = \bL(:, \sC_{j})$.  Since $\bL$ is lower
triangular, it follows that $\bL_{j}(1:c_j-1, :) = \bzero$.
We store the matrix $\bL_{j}(\sC_{j}, :)$ explicitly as a dense matrix
and refer to this as the supernode's {\em diagonal block}
$\bL^{D}_{j}$.  We also define $\sR_{j}$ to be the list of 
nonzero rows of $\bL_j$ below the diagonal block; that is,
\begin{equation}
  \sR_{j} = \left(
      k \in \sC^{O} :
        \exists p \in \sC_{j}, \, \ell_{k,p} \ne 0 \right)
      = \left( r_{j}^{1}, r_{j}^{2}, \ldots \right).
\end{equation}
%
Since columns in $\sC_{j}$ have similar non-zero patterns,
we store $\bL_{j}(\sR_{j}, :)$ as a dense matrix and refer to this as 
supernode $j$'s (compressed) off-diagonal block $\bL^{O}_{j}$.

Left-looking algorithms such as the one implemented in \cite{chen08} form
block columns $\bL_{j}$ in order from left to right.
We identify the {\em descendants} $\bbD_{j}$ of supernode $j$ as follows:
\begin{equation}
\label{eq:descendents}
  \bbD_{j} = \left\{
    1 \le k < j : \sR_{k} \cap \sC_{j} \ne \emptyset \right\};
\end{equation}
that is, supernodes from earlier in the factorization whose off-diagonal
row set intersects the column set of node $j$.  We also refer to node $j$
as an {\em ancestor} of node $k$ if $k \in \bbD_{j}$.
%
%
For convenience, we also define 
index lists relating rows in supernode $j$ to rows in a the
off-diagonal block $\bL^O_k$ of a descendant supernode $k$:
\begin{align}
  \rDkToj
    &= 
    \left( 1 \le p \le |\sR_{k}| : r^{p}_{k} \in \sC_{j} \right),
      \label{eq:rowsDiagKtoJ} \\
  \rOkToj
    &=
    \left( 1 \le p \le |\sR_{k}| : r^{p}_{k} \in \sR_{j} \right).
      \label{eq:rowsOffKtoJ}
\end{align}
Intuitively, \eqref{eq:rowsDiagKtoJ} 
helps us extract the rows of $\bL^O_k$
that influence the contents of $\bL^{D}_{j}$.
Similarly, 
we use \eqref{eq:rowsOffKtoJ} to extract rows of $\bL^O_k$
needed to form $\bL^O_j$.
We will also write $\srDkToj$
to denote the index list corresponding to $\rDkToj$ 
in the uncompressed structure,
$\left( r^{p}_{k} \in \sR_{k} : p \in \rDkToj \right)$, 
and similarly for $\srOkToj$.

We also find it convenient to define the function
$\scatterRows( \bB, R_{1}, R_{2} )$, where $R_{1} \subseteq R_{2}$ are ordered
index lists, and $\bB$ has $|R_{1}|$ rows.  This function returns a matrix
with $|R_{2}|$ rows by placing the contents of rows of $\bB$ in the output
according to the positions of entries of $R_{1}$ in $R_{2}$.  For example,
\begin{equation}
\label{eq:scatterEx}
  \scatterRows \left(
    \left( \begin{array}{cc} 1 & 2 \\ 3 & 4 \end{array} \right),
    \left\{ 3, 8 \right\},
    \left\{ 2, 3, 5, 8 \right\}
    \right)
  = \left( \begin{array}{cc}
    0 & 0 \\
    1 & 2 \\
    0 & 0 \\
    3 & 4
  \end{array} \right).
\end{equation}
We similarly define the functions $\scatterColumns( \bB, C_{1}, C_{2} )$, and
$\scatter( \bB, R_{1}, R_{2}, C_{1}, C_{2} )$ which composes
$\scatterRows$ and $\scatterColumns$.  Finally, we define $\gatherRows$
as a function which reverses the operation of $\scatterRows$; e.g., if
$\scatter(\bB, R_{1}, R_{2}) = \sC$ then $\gatherRows(\sC, R_{2}, R_{1}) = \bB$.

Forming the numerical
contents of supernode $j$ begins 
with assembly of a block column of the Schur complement:
\begin{align}
  \sU^{D}_{j} &= \bA(\sC_{j}, \sC_{j}) -
    \sum_{k \in \bbD_{j}}
      \scatter \left(
        \bL^{O}_{k} \left( \rDkToj, : \right) \!
        \bL^{O}_{k} \left( \rDkToj, : \right)^{T}
        \! \! , \,
        \srDkToj, \, \sC_{j}, \,
        \srDkToj, \, \sC_{j}
      \right)
      \label{eq:schurDiag} \\
  \sU^{O}_{j} &= \bA(\sR_{j}, \sC_{j}) -
    \sum_{k \in \bbD_{j}}
      \scatter \left(
        \bL^{O}_{k} \left( \rOkToj, : \right) \!
        \bL^{O}_{k} \left( \rDkToj, : \right)^{T}
        \! \! , \,
        \srOkToj, \, \sR_{j}, \,
        \srDkToj, \, \sC_{j}
      \right)
      \label{eq:schurOffDiag}
\end{align}
$\sU^{D}_{j}$ and $\sU^{O}_{j}$ are dense matrices with the same sizes
as $\bL^{D}_{j}$ and $\bL^{O}_{j}$.
We note that if node $k$ is a descendant of node $j$, then
$\sR_{k} \cap \sC^{O}_{j} \subseteq \sR_{j}$.
The Schur complement in node $j$ is formed by first extracting
dense row subsets of $\bL^{O}_{k}$ for each descendant $k$, then forming the
matrix products from (\ref{eq:schurDiag}-\ref{eq:schurOffDiag}) using dense
matrix arithmetic and finally scattering the result to $\sU^{D}_{j}$ and
$\sU^{O}_{j}$.  Next, the diagonal
and off-diagonal blocks of supernode $j$ are formed as follows:
\begin{align}
  \bL^{D}_{j} &= \chol \left( \sU^{D}_{j} \right) \\
  \bL^{O}_{j} &= \sU^{O}_{j} \left( \bL^{D}_{j} \right)^{-T}
\end{align}
This procedure is summarized in Algorithm~\ref{ALGfactorSupernode}.

\begin{algorithm}[!htb]

\SetKwData{Ii}{i}
\SetKwData{Dd}{d}
\SetKwData{Each}{each}
\SetKwData{DiagUpdate}{diagUpdate}
\SetKwData{OffDiagUpdate}{offDiagUpdate}
\SetKwInOut{Input}{input}
\SetKwInOut{Output}{output}

\Input{$\bA$, $j$, $\bL$, $\bbD_{j}$}
\Output{$\bLdj$, $\bLoj$}

\Begin{
  \SetAlgoVlined
  \CommentSty{// Initialize the Schur complement blocks} \\[-2pt]
  $\sUdj \longleftarrow \bA(\sC_{j}, \sC_{j})$ \\
  $\sUdj \longleftarrow \bA(\sR_{j}, \sC_{j})$ \\[5pt]
  
  \For{$\Each \, k \in \mathbb{D}_{j}$}{
    \CommentSty{// Build dense update blocks} \\[-2pt]
    $\DiagUpdate \longleftarrow \bLok(\rDkToj,:) * \bLok(\rDkToj,:)^{T}$ \\
    $\OffDiagUpdate \longleftarrow \bLok(\rOkToj,:) * \bLok(\rDkToj,:)^{T}$ \\
      [3pt]
    \CommentSty{// Scatter updates to the Schur complement} \\[-2pt]
    $\sUdj \longleftarrow \sUdj
      - \scatter(\DiagUpdate, \srDkToj, \sC_{j}, \srDkToj, \sC_{j})$ \\
    $\sUoj \longleftarrow \sUoj
      - \scatter(\OffDiagUpdate, \srOkToj, \sR_{j}, \srDkToj, \sC_{j})$ \\
  }
  \CommentSty{// Factor node $j$'s diagonal block} \\[-2pt]
  $\bLdj \longleftarrow \cholesky(\sUdj)$ \\[5pt]
  \CommentSty{// Dense triangular solve} \\[-2pt]
  $\bLoj \longleftarrow \sUoj (\bLdj)^{-T}$ \\[5pt]

  $\Return \, \, \, \bLdj, \bLoj$
}

\cprotect\caption{{$\factorSupernode$:
  Computes the diagonal and off-diagonal blocks $\bL^{D}_{j}$ and
  $\bL^{O}_{j}$ for supernode $j$.  Provided inputs are the matrix to
  be factored, as well as the partially constructed factor $\bL$.
  every supernode $k \in \bbD_{j}$ (node $j$'s descendants) is assumed
  to have already been factored.
  }}
\label{ALGfactorSupernode}
\end{algorithm}

\subsection{Fill-Reducing Ordering}
\label{SECordering}

\begin{figure}
  \begin{center}
  \includegraphics[width=1.45in]{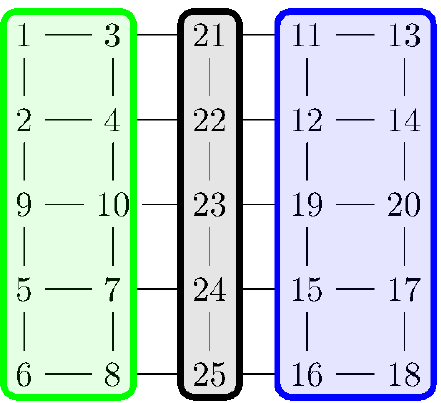}
  \hspace{1.5cm}
  \includegraphics[width=2.9in]{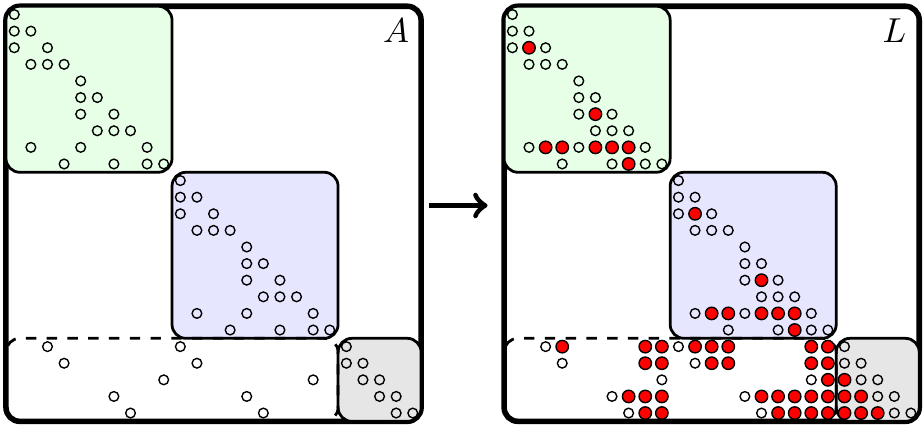}
  \end{center}
  \caption{Fill during nested-dissection Cholesky factorization for a 2D
           mesh example (left).  When nodes in a mesh are ordered
           so that a vertex separator appears last after the subdomains
           it separates (right), fill in the Cholesky factor (indicated by
           red circles) is restricted
           to the diagonal blocks corresponding to interactions within each
           subdomain and within the separator and to off-diagonal blocks
           associated with subdomain-separator interactions.}
  \label{fig:nestedDissection}
\end{figure}

Figure \ref{fig:nestedDissection} provides an example of how a nested
dissection ordering might be used on a simple two-dimensional domain, and
how this ordering influences fill in the Cholesky factor of the reordered
matrix.

\section{Off-Diagonal Block Compression}
\label{SECoffdiagcomp}

In \S\ref{SECsupernodal} we introduced the notion of supernodes with
dense diagonal and off-diagonal blocks $\bL^{D}_{j}$ and $\bL^{O}_{j}$.
In this section, we discuss the process of approximating off-diagonal
blocks $\bL^{O}_{j}$ with low-rank matrices.  The matrix $\bLoj$ stores
interactions between supernode $j$, and other supernodes in occurring later
in the factorization.

\subsection{Block Selection and Ordering}
\label{SECoffdiagordering}

We use nested dissection~\cite{george73} to construct a fill-reducing ordering
as originally discussed in \S\ref{SECordering}.  We choose nested dissection
because the geometric structure introduced by this method yields a factor
matrix $\bL$ in which many dense submatrices are amenable to low-rank
approximation.  We will discuss this property in more detail in
\S\ref{SECodlowrank}.
Henceforth, we will
assume that the matrix $\bA$ has already been symmetrically permuted using
nested dissection.
We expect that the interactions between large separators in the factorization
will exhibit rapidly decaying rank structure (see \S\ref{SECodlowrank}).
Therefore, we represent each ``large'' separator
from the nested dissection hierarchy with a supernode.  The off-diagonal
blocks $\bLoj$ associated with these supernodes describe interactions between
large separators in the factor (see Figure~\ref{FIGoffdiagBlocks}).

Standard sparse Cholesky solvers such as CHOLMOD \cite{chen08} may optionally
use nested
dissection for reordering.  However, the process of constructing supernodes
used by these solvers does not guarantee a one-to-one relationship between
supernodes and separators from the nested dissection ordering.  In our
algorithm, we choose a tolerance $\tau_{O} \in \bbZ^{+}$ and introduce a
supernode for every separator with at least $\tau_{O}$ variables.  The remaining
indices in our reordered matrix $\bA$ are gathered in to supernodes using
methods identical to \cite{chen08}.

\begin{figure}
  \begin{center}
  \includegraphics[width=\hsize]{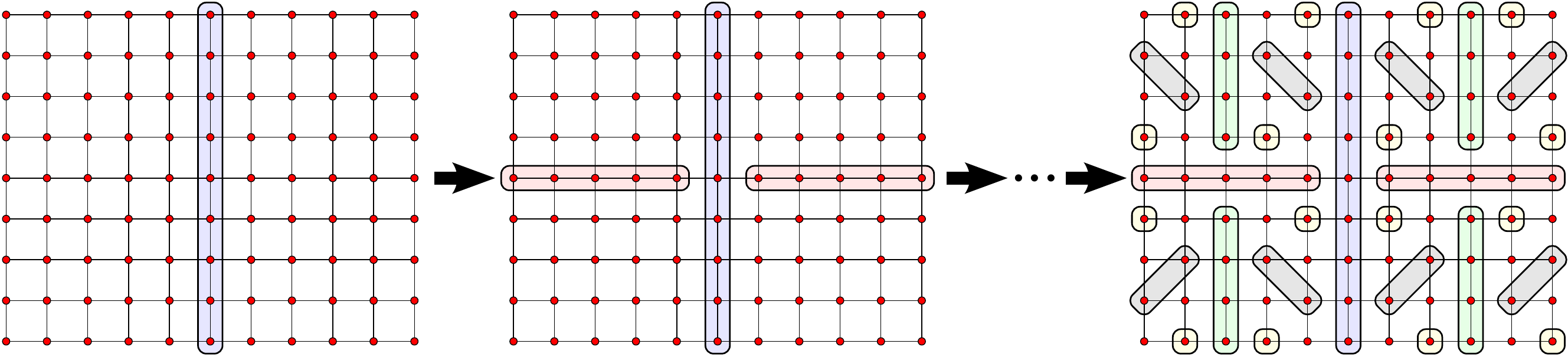}
  \end{center}
  \caption{Nested dissection is applied to a regular, two-dimensional grid.
           At each level, the domain is recursively subdivided by the
           introduction of separators.  Each image in the sequence depicts a
           level in the recursive dissection of the domain.
           Separators at the same level in the
           nested dissection hierarchy are rendered with the same color.  In the
           rightmost figure we see that five levels of nested dissection
           fully decompose this domain in to subdomains of unit size.}
  \label{FIGnestedSequence}
\end{figure}

\subsection{Block Compression Scheme}
\label{SECodblockcomp}

Consider the state of the factorization
immediately prior to forming the factor contents for supernode $j$:
\begin{equation}
\label{eq:preNodeFactor}
  \left( \begin{array}{ccc}
    \bA^{D}_{pre} & \multicolumn{2}{c}{sym} \\
    \multirow{2}{*}{$\bA^{O}_{pre}$} & \bA^{D}_{j} & sym \\
    & \bA^{O}_{j} & \bA_{post}
    \end{array} \right)
  =
  \left( \begin{array}{ccc}
    \bL^{D}_{pre} & & \\
    \multirow{2}{*}{$\bL^{O}_{pre}$} & \bI & \\
    & \bzero & \bI
    \end{array} \right)
  \left( \begin{array}{ccc}
    \bL^{D}_{pre} & & \\
    \multirow{2}{*}{$\bL^{O}_{pre}$} & \sU^{D}_{j} & \\
    & \sU^{O}_{j} & \sU_{post}
    \end{array} \right)^{T}
\end{equation}
Here {\em pre} and {\em post} refer to the sets of columns occurring before
and after supernode $j$, respectively.  Note that the Schur complement
$\sU_{post}$ is never formed explicitly since we only form Schur complements
one supernode at a time.  We also note that given the definition of
$\sU_{j}$ and $\bL_{j}$, it is necessary to apply the $\scatter$ operator to
these matrices to make \eqref{eq:preNodeFactor} valid, but this has been
omitted here for brevity.
Following factorization of node $j$, we have:
\begin{equation}
  \label{eq:postNodeFactor}
   \left( \begin{array}{ccc}
    \bA^{D}_{pre} & sym & sym \\
    \multirow{2}{*}{$\bA^{O}_{pre}$} & \bA^{D}_{j} & sym \\
    & \bA^{O}_{j} & \bA_{post}
    \end{array} \right)
  =
  \left( \begin{array}{ccc}
    \bL^{D}_{pre} & & \\
    \multirow{2}{*}{$\bL^{O}_{pre}$} & \bL^{D}_{j} & \\
    & \bL^{O}_{j} & \bI
    \end{array} \right)
  \left( \begin{array}{ccc}
    \bL^{D}_{pre} & & \\
    \multirow{2}{*}{$\bL^{O}_{pre}$} & \bL^{D}_{j} & \\
    & \bL^{O}_{j} & \widetilde{\sU}_{post}
    \end{array} \right)^{T} 
\end{equation}
where
\begin{equation}
\label{eq:scOriginal}
  \widetilde{\sU}_{post} = \sU_{post}
    - \bL^{O}_{j} \left( \bL^{O}_{j} \right)^{T}.
\end{equation}
Assuming 
$\bA$ is positive definite, 
$\widetilde{\sU}_{post}$ must also be positive definite.

If we approximate the off-diagonal part of supernode $j$
with a low-rank matrix -- ${\bLoj \approx \bV \bU^{T}}$ -- then the
Schur complement in \eqref{eq:scOriginal} is approximated by
\begin{equation}
\label{eq:scModified}
  \overline{\sU}_{post} = \sU_{post}
    - \bV \bU^{T} \bU \bV^{T}.
\end{equation}
We choose $\bV$ and $\bU$ using a method similar to \cite{Xia2012c} so that
this modified Schur complement \eqref{eq:scModified} is guaranteed to remain
positive definite.
Namely, we choose $\bU$
to have orthonormal columns and $\bV$ to be the projection of $\bL^{O}_{j}$
on to this basis; ${\bV = \bLoj \bU}$.  We can write
$\bLoj = [ \bV \, \, \, \overline{\bV} ] [ \bU \, \, \, \overline{\bU} ]^{T}$
where $\overline{\bU}$ is a (non-unique) matrix with orthonormal columns,
$\bU^{T} \, \overline{\bU} = \bzero$, and
$\overline{\bV} = \bLoj \overline{\bU}$.
Given these properties, we can rewrite \eqref{eq:scOriginal} as
\begin{equation}
  \widetilde{\sU}_{post} =
    \sU_{post} - \bV \bV^{T} - \overline{\bV} \, \overline{\bV}^{T}
\end{equation}
and \eqref{eq:scModified} as
\begin{align}
  \overline{\sU}_{post}
    &=
    \sU_{post} - \bV \bV^{T} \\
    &=
    \widetilde{\sU}_{post} + \overline{\bV} \, \overline{\bV}^{T}.
      \label{eq:scPlusPos}
\end{align}
Since $\widetilde{\sU}_{post}$ is positive definite and
$\overline{\bV} \, \overline{\bV}^{T}$ is positive semi-definite, it
follows from \eqref{eq:scPlusPos} that $\overline{\sU}_{post}$ remains
positive definite under this approximation.

\subsection{Low-Rank Structure}
\label{SECodlowrank}

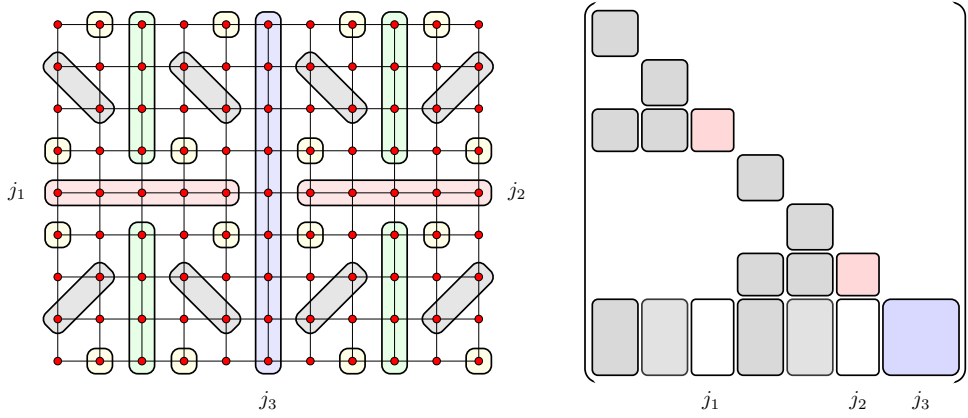
\begin{figure}
  \resizebox{\hsize}{!}{%
  \begin{tikzpicture}
  \node [anchor=south west] (label) at (0,0)
    {\includegraphics{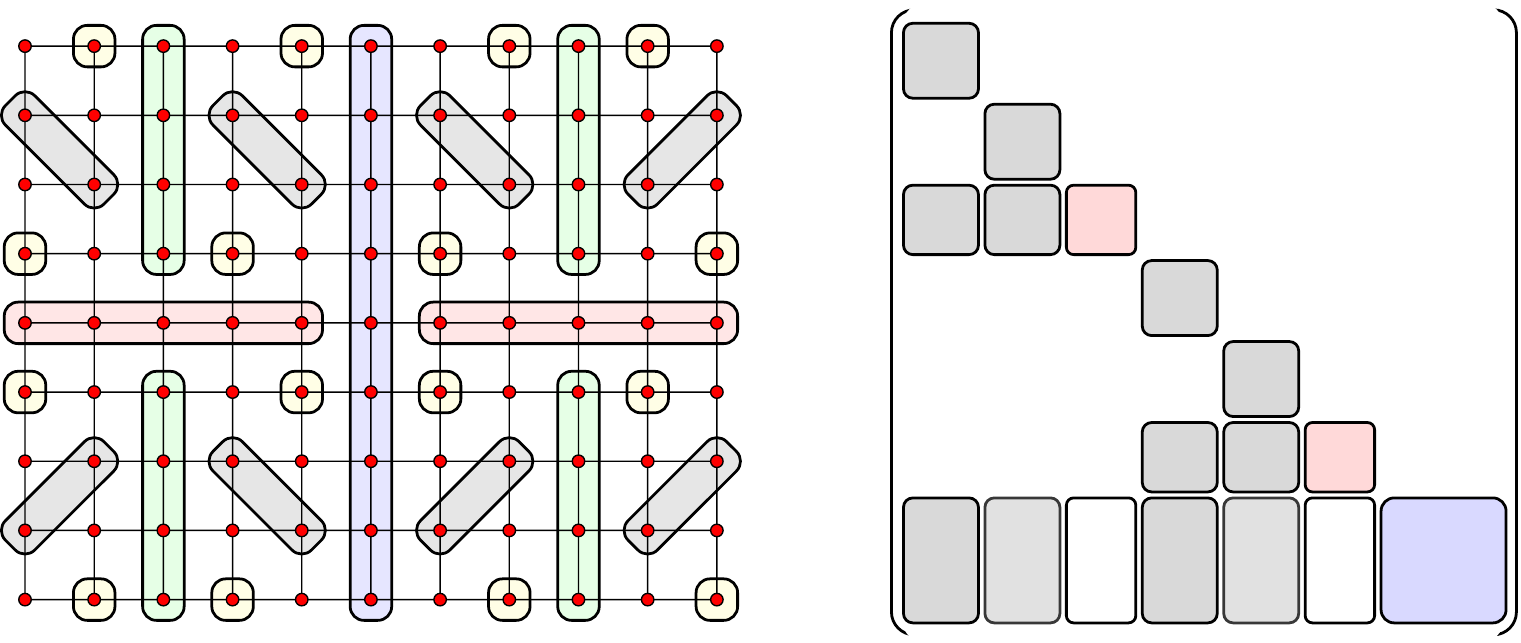}};
  \draw (-0.3,3.4) node {$j_{1}$};
  \draw (8.05,3.4) node {$j_{2}$};
  \draw (3.9, -0.1) node {$j_{3}$};
  \draw (11.29, -0.1) node {$j_{1}$};
  \draw (13.75, -0.1) node {$j_{2}$};
  \draw (14.80, -0.1) node {$j_{3}$};
  \end{tikzpicture}
  }%
  \caption{Consider the two-dimensional grid on the left, decomposed via nested
           dissection (see Figure~\ref{FIGnestedSequence}).  If we set
           $\tau_{O} = 5$, then the three largest separators will be identified
           as supernodes (labelled $j_{1}$, $j_{2}$ and $j_{3}$).  The resulting
           matrix structure is shown on the right.  The off-diagonal blocks
           $\bL^{O}_{j_{1}}$ and $\bL^{O}_{j_{2}}$ are shown in white and
           describe interactions between separator $j_{3}$ and separators
           $j_{1}$ and $j_{2}$, respectively.  In the domain picture on the
           left we see that, due to the way these separators intersect
           geometrically, these interactions tend to mostly occur over large
           distances.  This justifies the use of low-rank matrices to
           approximate these interactions.  The parts of the factor for which
           we do not apply any compression are shown in gray.  These blocks
           can be evaluated either using the standard supernodal factorization
           algorithm, or using the interior blocks approach discussed in
           \S\ref{SECinteriorBlocks}.
           }
  \label{FIGoffdiagBlocks}
\end{figure}

Recall from \S\ref{SECoffdiagordering} that we 
use nested dissection
to reduce fill and represent large separators
as supernodes in the factorization.
Nested dissection orderings are built
entirely based on $\bA$'s graph structure.
However, in problems defined on physical domains (say, discretizations of partial
differential equations on two- or three-dimensional domains) separators
also have a convenient geometric interpretation.
In these problems, separators are geometric regions which bisect subdomains
of the original problem domain (see Figures~\ref{FIGnestedSequence} and
\ref{FIGoffdiagBlocks}).  For example, in many three-dimensional problems the
nested dissection separators are surfaces which cut the domain in to disjoint
pieces.
Given a partial factorization of a matrix $\bA$, the remaining Schur
complement $\sU$ behaves like a discretization of a boundary
integral equation \cite{chandrasekaran2010}.
For many problems, these discretizations
will yield smooth coefficients for matrix indices which are geometrically
distant from each other in the original problem domain.
The structure of separators produced by nested dissection tends to ensure
that the interactions between pairs of large separators occur mostly over
large distances, with only a handful of ``near-field'' interactions
(see Figure~\ref{FIGoffdiagBlocks}).  If supernode $j$ corresponds to a large
separator and $\sUoj$ is this node's off-diagonal block in the Schur complement,
then we expect that $\sUoj$ should have rapidly decaying rank structure due
to the property discussed above.  As such, $\sUoj$ (and, likewise, $\bLoj$)
admits a low-rank approximation.
Therefore, we compress the off-diagonal blocks in supernodes/separators which
are sufficiently large (larger than $\tau_{O}$).  We use the following
notation for this low-rank approximation:
\begin{equation}
  \bLoj \approx \bV_{j} \bU_{j}^{T} \quad \textrm{where} \quad
  \bV_{j} \in \bbR^{|\sR_{j}| \times q},
  \bU_{j} \in \bbR^{|\sC_{j}| \times q}  \quad \textrm{and} \quad
  q \ll |\sR_{j}|, |\sC_{j}|
\end{equation}

%
%

\subsection{Compression Algorithm}
\label{SECodcompalg}

\begin{algorithm}[!htb]

\SetKwData{Ii}{i}
\SetKwData{Dd}{d}
\SetKwData{Each}{each}
\SetKwData{DiagUpdate}{diagUpdate}
\SetKwData{OffDiagUpdate}{offDiagUpdate}
\SetKwData{Transpose}{transpose}
\SetKwData{Boolean}{boolean}
\SetKwData{True}{true}
\SetKwData{False}{false}
\SetKwFunction{Cholesky}{Cholesky}
\SetKwInOut{Input}{input}
\SetKwInOut{Output}{output}

\Input{$\bA$, $j$, $\bL$, $\bbD_{j}$, $\bG$}
\Output{$\bB = \bLoj \bG$}

\Begin{
  \SetAlgoVlined
  \CommentSty{// Apply node $j$'s diagonal inverse to the input} \\[-2pt]
  $\bG \longleftarrow \diagonalSolve(j, \bG, \Transpose = \True)$ \\[5pt]

  \CommentSty{// Multiply by the desired block from $\bA$}
    \\[-2pt]
  $\bW \longleftarrow \bA(\sR_{j}, \sC_{j}) \bG$ \\[5pt]
  
  \For{$\Each \, k \in \mathbb{D}_{j}$}{
    \CommentSty{// Extract the required sub-matrix from $\bG$} \\[-2pt]
    $\bG_{sub} \longleftarrow \gatherRows(\bG, \sC_{j}, \srDkToj)$ \\[5pt]

    \CommentSty{// Form the needed product with two multiplications} \\[-2pt]
    $\bT \longleftarrow \left[ \bLok(\rDkToj, :) \right]^{T} \bG_{sub}$ \\
    $\bT \longleftarrow \bLok(\rOkToj, :) \bT$ \\[5pt]

    \CommentSty{// Scatter result to the output matrix} \\[-2pt]
    $\bW \longleftarrow \bW - \scatterRows(\bT, \srOkToj, \sR_{j})$
  }

  $\Return \, \, \, \bW$
}

\cprotect\caption{{$\offDiagonalMultiply$:
  Computes the product $\bB = \bLoj \bG$ given some input matrix
  $\bG \in \bbR^{|\sC_{j}| \times r}, r > 0$.
  The function $\diagonalSolve(j, \bX, \Transpose)$
  applies the inverse of $\bLdj$ to the input matrix.  If transpose is
  set to $\True$, then $\diagonalSolve$ forms the product
  $(\bLdj)^{-T} \bX$ instead.  See \S\ref{SECdiagcomp} and Algorithm
  \ref{ALGdiagonalSolve} for a detailed description of $\diagonalSolve$.
  }}
\label{ALGoffDiagMult}
\end{algorithm}

Next we discuss our approach for forming the low-rank approximation
$\bLoj \approx \bVj \bUj^{T}$.
Our compression strategy must satisfy two requirements:
\begin{enumerate}
  \item The off-diagonal block $\bLoj$ may be expensive to construct and store.
        Therefore, we wish
        to build $\bVj$, $\bUj$
        without explicitly constructing 
        $\bLoj$.

  \item The original supernode factorization procedure presented in algorithm
        \ref{ALGfactorSupernode} makes effective use of dense matrix
        arithmetic, allowing for very efficient
        implementations~\cite{chen08}.
        Our compression algorithm should preserve this property.
\end{enumerate}

To satisfy these requirements, we 
use
randomized low-rank approximation
algorithms~\cite{liberty07,Halko2011}.  Similar randomized methods have been
used previously 
in
rank-structured sparse solvers~\cite{Engquist2010,Xia2012b}.
The key insight behind these randomized
algorithms is that a ``good'' rank-$q$ approximation to a matrix $\bB$ can
be found by considering products of the form $\bB \bG$, where $\bG$ is a
randomly generated matrix with $q + p$ columns and $p > 0$ is a small
oversampling parameter (typically $p \approx 5-10$ is suitable).  In this
case, we consider a rank-$q$ approximation to be good if it is close (in the
2-norm) to the best rank-$q$ 
approximation provided by $\bB$'s singular value decomposition (SVD).  
If the singular values of $\bB$ decay slowly, then
obtaining such an approximation may require
us to instead form products of the
form $\bC = (\bB \bB^{T})^{s} \bB \bG$, where $s \ge 0$ is a small number of
power iterations.
The theory behind these methods states that $\bC$ provides
a column basis for a low-rank approximation of $\bB$ which is close to optimal.
Moreover, constructing $\bC$ only requires a small number of
matrix multiplications involving $\bB$ and $\bB^{T}$.  Therefore, requirement
1 above is satisfied.  We can also perform these multiplications in a way
that leverages dense matrix arithmetic similar to
Algorithm~\ref{ALGfactorSupernode}, satisfying requirement 2.

Algorithm~\ref{ALGoffDiagMult} efficiently forums products $\bLoj \bG$ for
an arbitrary dense matrix $\bG$.  As discussed earlier, we also require
products of the form $\left( \bLoj \right)^{T} \bG$.  These products are
formed by the function $\offDiagonalMultiplyTrans$.  This function has
similar structure fo Algorithm~\ref{ALGoffDiagMult}.
Finally, Algorithm \ref{ALGapproxOffDiag} uses
the $\offDiagonalMultiply$ and $\offDiagonalMultiplyTrans$ functions
to form a low-rank approximation
$\bLoj \approx \bVj \bUj^{T}$ for node $j$'s off-diagonal block,
where $\bUj$ is chosen to have orthonormal columns (as discussed in
\S\ref{SECodblockcomp}).
We note that Algorithm~\ref{ALGoffDiagMult}
assumes the matrix $\bLok$ is stored explicitly for
all descendants $k \in \bbD_{j}$.  In practice, some of these blocks
may also have been assigned low-rank representations
$\bLok \approx \bVk \bUk^{T}$.  If this is the case, then we replace lines
10-11 in Algorithm~\ref{ALGoffDiagMult} with
\begin{align}
  \bU_{prod}  & \longleftarrow \bUk^{T} \bUk \notag \\
  \bT         & \longleftarrow \left[ \bVk(\rDkToj, :) \right]^{T} \bG_{sub}
    \notag \\
  \bT         & \longleftarrow \bU_{prod} \bT \notag \\
  \bT         & \longleftarrow \bVk(\rOkToj, :) \bT.
    \label{eq:lowRankMult}
\end{align}

\begin{algorithm}[!htb]

\SetKwData{Ii}{i}
\SetKwData{Dd}{d}
\SetKwData{Each}{each}
\SetKwData{DiagUpdate}{diagUpdate}
\SetKwData{OffDiagUpdate}{offDiagUpdate}
\SetKwData{Transpose}{transpose}
\SetKwData{Boolean}{boolean}
\SetKwData{True}{true}
\SetKwData{False}{false}
\SetKwFunction{Cholesky}{Cholesky}
\SetKwInOut{Input}{input}
\SetKwInOut{Output}{output}

\Input{$\bA$, $j$, $\bL$, $\bbD_{j}$,
       off-diagonal rank $s_{j}$, number of power iterations $q$}
\Output{$\bVj, \bUj$ such that $\bLoj \approx \bVj \bUj^{T}$}

\Begin{
  \SetAlgoVlined
  \CommentSty{// Initialize a random matrix} \\[-2pt]
  $\bG \longleftarrow \randomMatrix(|\sR_{j}|, s_{j})$ \\[5pt]

  \CommentSty{// Implicitly form the product $\bLoj \bG$} \\[-2pt]
  $\bG \longleftarrow \offDiagonalMultiplyTrans(\bA, j, \bL, \bbD_{j}, \bG)$
    \\[5pt]
  
  \CommentSty{// Run additional power iterations} \\[-2pt]
  \For{$i = 1$ to $q$}{
    $\bG \longleftarrow \offDiagonalMultiply(\bA, j, \bL, \bbD_{j}, \bG)$ \\
    $\bG \longleftarrow \offDiagonalMultiplyTrans(\bA, j, \bL, \bbD_{j}, \bG)$ \\
  }

  \CommentSty{// Extract an orthonormal row basis for $\bLoj$} \\[-2pt]
  $\bU_{j} \longleftarrow \makeOrthonormal(\bG)$ \\[5pt]

  \CommentSty{// Compute $\bVj$ by projecting on to $\bUj$} \\[-2pt]
  $\bVj \longleftarrow \offDiagonalMultiply(\bA, j, \bL, \bbD_{j}, \bUj)$ \\[5pt]

  $\Return \, \, \, \bVj, \bUj$
}

\cprotect\caption{{$\approximateOffDiagonal$:
  Builds a low-rank approximation $\bLoj \approx \bVj \bUj^{T}$ for node
  $j$'s off-diagonal block.  This algorithm assumes that node $j$'s diagonal
  block has already been factored.  The function $\randomMatrix(m, n)$
  generates an $m \times n$ matrix whose entries are drawn from a Gaussian
  distribution with mean 0 and unit variance.  The function $\makeOrthonormal$
  returns an orthonormal basis for the column space of its input matrix.
  This can be accomplished by means of -- for example -- a QR factorization.}}
\label{ALGapproxOffDiag}
\end{algorithm}

\section{Diagonal Block Compression}
\label{SECdiagcomp}

Section \S\ref{SECoffdiagcomp} discussed the process of constructing
a sparse Cholesky factorization in which off-diagonal interactions between
large separators are approximated with low-rank matrices.
In large, three-dimensional problems, larger separators may include thousands
to tens of thousands of variables.
For these problems, the compression scheme
from \S\ref{SECoffdiagcomp} can provide a significant reduction in both
memory usage over standard factorizations, while still providing a factor
that serves as an excellent preconditioner.
However, if supernode $j$ is large,
then building the dense diagonal matrix $\bL^{D}_{j}$ may also require
significant computation and storage.  In this section, we discuss an 
approach to compressing diagonal blocks $\bL^{D}_{j}$.

\subsection{Low-Rank Structure}
\label{SECexposelowrank}

In \S\ref{SECoffdiagcomp} we saw that low-rank behavior in
off-diagonal blocks $\bLoj$ is exposed by the geometric structure
of nested dissection.  We can reorder variables within a separator to expose
similar low-rank structure within diagonal blocks $\bLdj$. 
Since all columns in a given supernode are treated as having the same fill
pattern, we are can perform this reordering without affecting
accuracy, memory usage, or computation time.
Consider the top-level separator shown in figure
\ref{FIGnestedSequence}.
Suppose that
the indices of vertices in this separator are ordered sequentially from top
to bottom and that the $9 \times 9$ diagonal block for this separator is
written as a $2 \times 2$ block matrix (assume, without loss of generality,
that the first block row/column has four entries, and that the second has five):
\begin{equation*}
  \bLdj = \left( \begin{array}{cc}
    \bL_{11} & \bzero \\
    \bL_{21} & \bL_{22}
  \end{array} \right)
\end{equation*}
As a result of the ordering discussed above, $\bL_{21}$ stores interactions
between vertices
in the top half of the separator with vertices in the bottom half.
As we discussed in \S\ref{SECodlowrank},
the spatial separation between these groups of variables suggests that $\bL_{21}$
can be approximated with a low-rank matrix $\bL_{21} \approx \bV \bU^{T}$.
We can apply this argument recursively to $\bL_{11}$ and $\bL_{22}$
to achieve further compression.

The example above assumes that the rows/columns of $\bLdj$ are
ordered 
such that off-diagonal blocks of $\bLdj$ exhibit low-rank structure.  
Finding such an ordering
is straightforward when $\bA$ 
comes from a PDE discretization on a regular mesh
like the one pictured in Figure~\ref{FIGnestedSequence}.
However, obtaining a suitable ordering for general,
three-dimensional problems on irregular domains is nontrivial.

\begin{figure}
  \resizebox{\hsize}{!}{%
  \begin{tikzpicture}
  \node [anchor=south west] (label) at (0,0)
    {
      \includegraphics[width=1.8in]{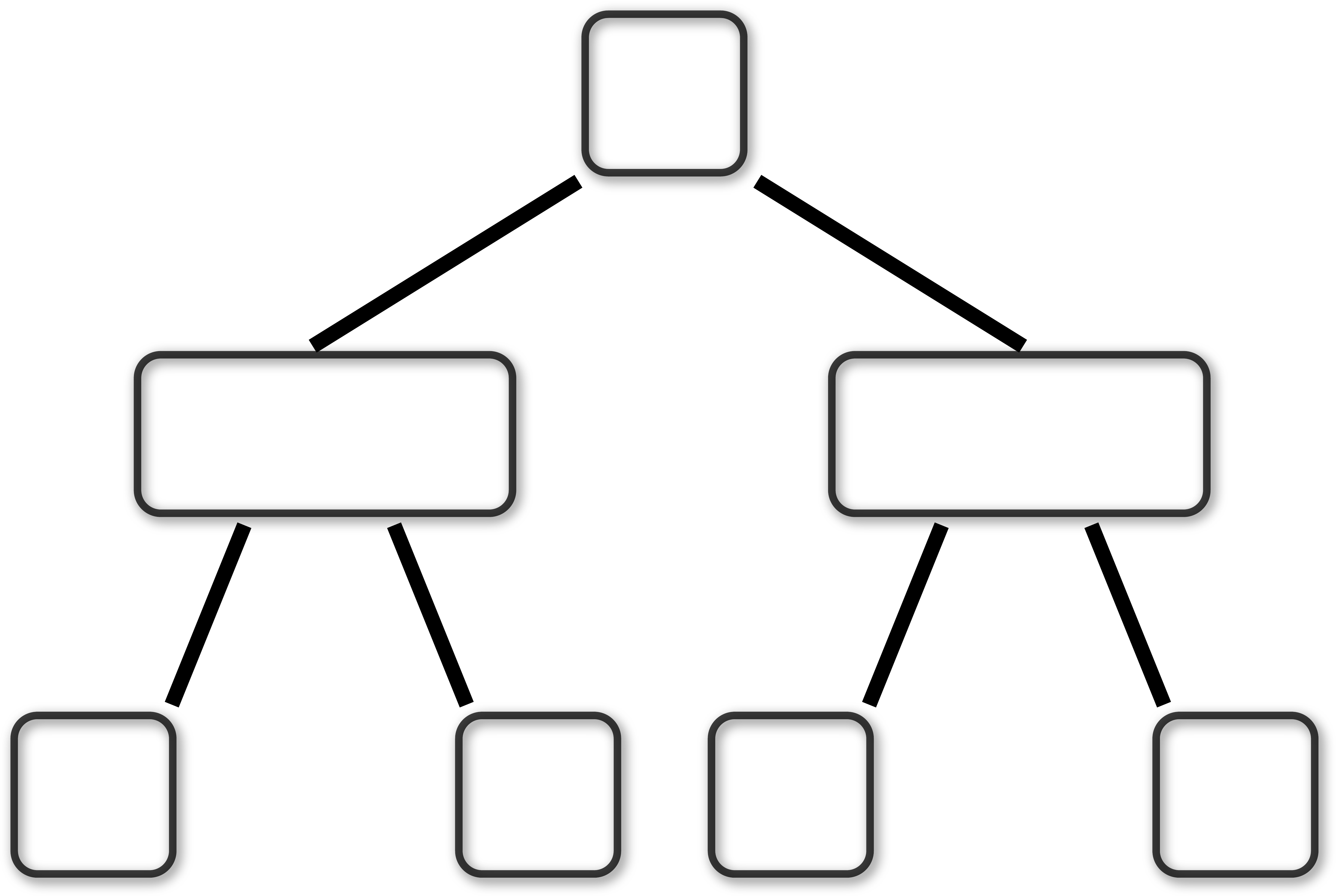}
    };

  \draw (10.0,1.7) node {$
  \begingroup \renewcommand*{\arraystretch}{1.5}
    \bAdj \longrightarrow
    \left( \begin{array}{cc:cc}
      \multicolumn{1}{c:}{\bAdj(I_{1}, I_{1})} &
        \multicolumn{1}{c:}{sym.} & 
        \multicolumn{2}{c}{\multirow{2}{*}{sym.}} \\ \cdashline{1-2}
      \multicolumn{1}{c:}{\bAdj(I_{2}, I_{1})} &
        \multicolumn{1}{c:}{\bAdj(I_{2}, I_{2})} &  &  \\ \hdashline
      \multicolumn{2}{c:}{
        \multirow{2}{*}{
          $\bAdj(I_{3} \cup I_{4}, I_{1} \cup I_{2})$}} &
            \multicolumn{1}{c}{\bAdj(I_{3}, I_{3})}
            & \multicolumn{1}{:c}{sym.} \\ \cdashline{3-4}
      & & \multicolumn{1}{c}{\bAdj(I_{4}, I_{3})} &
        \multicolumn{1}{:c}{\bAdj(I_{4}, I_{4})}
    \end{array} \right)
  \endgroup
  $};

  \draw (2.4,2.85) node {$\sC_{j}$};
  \draw (1.23,1.7) node {$I_{1} \cap I_{2}$};
  \draw (3.63,1.7) node {$I_{3} \cap I_{4}$};
  \draw (0.45,0.45) node {$I_{1}$};
  \draw (1.95,0.44) node {$I_{2}$};
  \draw (2.83,0.45) node {$I_{3}$};
  \draw (4.35,0.45) node {$I_{4}$};
  \end{tikzpicture}
  }
  \caption{We recursively partition supernode $j$'s variable indices $\sC_{j}$
           resulting in the tree structure shown on the left (assuming two
           levels of partitioning in this case).
           We permute the indices in $\sC_{j}$ so that, following this
           permutation, the index blocks associated with leaves in this tree
           structure appear sequentially along the diagonal of $\bAdj$.
           The resulting permutation to $\bAdj$ for this two-level example is
           shown on the right.
           }
  \label{FIGdiagOrdering}
\end{figure}

Recall from \S\ref{SECodlowrank} that separators in the nested dissection
hierarchy are geometric regions partition the original problem domain.
In three dimensions, we intuitively expect these separators
to look like two-dimensional surfaces inside of the original problem domain.
Our solver exposes low-rank structure in $\bLdj$ by reordering indices within
separators 
so that off-diagonal blocks in $\bLdj$ describe interactions between
spatially separated pieces of the separator region.
We begin by assigning a three-dimensional position $\bx_{i} : i \in \sC_{j}$
to each index associated with supernode/separator $j$.
There are many techniques for spatially partitioning the positions $\bx_{i}$.
Currently, we use a simple axis-based splitting scheme.
We partition the positions $\bx_{i}$ in to two subsets by sorting them along
the longest bounding box axis of the set
$\{\bx_{i} : i \in \sC_{j} \}$ and splitting this sorted list in to two
equal-sized pieces.
This process is applied recursively until the separator has been partitioned
in to subdomains with at most $\tau_{D}$ variables.
The paramter $\tau_{D}$ is chosen in advance
as the size of the largest diagonal block that we wish to represent explicitly
in the factor matrix.
We use this partitioning to reorder the indices within a supernode in a way
that exposes low-rank structure.
See Figure~\ref{FIGdiagOrdering} for an explanation of how this permutation
is built.  Using the two-level partitioning example shown this figure,
we label blocks of the diagonal factor block $\bLdj$ as follows:
\begin{equation}
\begingroup \renewcommand*{\arraystretch}{1.5}
  \bLdj =
  \left( \begin{array}{cc:cc}
    \multicolumn{1}{c:}{\dBlock{1}} &
      \multicolumn{1}{c:}{\bzero} &
        \multicolumn{2}{c}{\multirow{2}{*}{$\bzero$}} \\ \cdashline{1-2}
  \multicolumn{1}{c:}{\dBlock{2}} & \dBlock{3} & & \\ \hdashline
  \multicolumn{2}{c:}{\multirow{2}{*}{$\dBlock{4}$}} &
    \multicolumn{1}{c}{\dBlock{5}} &
      \multicolumn{1}{:c}{\bzero} \\ \cdashline{3-4}
  & & \dBlock{6} & \multicolumn{1}{:c}{\dBlock{7}}
  \end{array} \right).
  \label{eq:ldjLabel}
\endgroup
\end{equation}
Blocks are numbered in the order in which they must be formed
during factorization (intuitively, top to bottom and left to right).
In this $4 \times 4$
example, $\bLdj$ is be approximated as follows:
\begin{equation}
\begingroup \renewcommand*{\arraystretch}{1.5}
  \bLdj \approx
  \left( \begin{array}{cc:cc}
  \multicolumn{1}{c:}{\dBlock{1}} &
    \multicolumn{1}{c:}{\bzero} &
      \multicolumn{2}{c}{\multirow{2}{*}{$\bzero$}} \\ \cdashline{1-2}
  \multicolumn{1}{c:}{\vBlock{2} (\uBlock{2})^{T}} &
    \multicolumn{1}{c:}{\dBlock{3}} & & \\ \hdashline
  \multicolumn{2}{c:}{\multirow{2}{*}{$\vBlock{4} (\uBlock{4})^{T}$}}
    & \dBlock{5} &
      \multicolumn{1}{:c}{\bzero} \\ \cdashline{3-4}
  & & \vBlock{6} (\uBlock{6})^{T} &
    \multicolumn{1}{:c}{\dBlock{7}}
  \end{array} \right).
  \label{eq:ldjComp}
\endgroup
\end{equation}
For block
$s$ in this matrix, let $\rDs$ and $\cDs$ be the set of rows and columns over
which block $s$ is defined, relative to $\bLdj$.
That is, $\bLdj(\rD{s}, \cD{s}) = \dBlock{s}$.
Similarly, let $\srDs$ and $\scDs$
refer to the same row and column sets, but relative to the entire factor
$\bL$, so that $\bL(\srDs, \scDs) = \dBlock{s}$.

In summary, we permute the original
matrix $\bA$ in two main steps.  The first is fill-reducing ordering using
nested dissection.  As we noted in \S\ref{SECoffdiagcomp},
this step exposes low-rank structure in certain off-diagonal submatrices
of $\bL$, allowing for compression.  The second stage of this permutation
consists of reordering indices {\em within} certain supernodes formed in the
first stage -- namly, those associated with large separators.  This does not
alter the sparsity of $\bL$, but does allow us to compress certain off-diagonal
submatrices of these large diagonal blocks.

\subsection{Compression Algorithm}
\label{SECdiagcompalg}

Next, we 
turn to compression of
off-diagonal blocks within a diagonal matrix $\bLdj$.
As in the compression methods discussed in \S\ref{SECoffdiagcomp}, we 
use randomized methods to construct
low-rank matrix approximations.
First, we define the $\diagonalSolve$ function, originally introduced in
Algorithm~\ref{ALGoffDiagMult}.
We consider a slight variation on this function, in which the
index $s$ of a block from $\bLdj$ is provided as an argument.
Invoking this function with argument $s$ solves a system of equations using the
smallest diagonal sub-block of $\bLdj$ containing $\dBlock{s}$.
Using
\eqref{eq:ldjComp} as an example, calling $\diagonalSolve$ with
$s = 2$ would solve a system using the inverse of
\begin{equation}
\begingroup \renewcommand*{\arraystretch}{1.5}
  \left( \begin{array}{cc}
    \dBlock{1} & \bzero \\ \vBlock{2} \left( \uBlock{2} \right)^{T} & \dBlock{3}
  \end{array} \right).
\endgroup
\end{equation}
For brevity, 
calling
$\diagonalSolve$ with no ``$s$'' argument
(as in Algorithm~\ref{ALGoffDiagMult}) solves
systems using the entire matrix $\bLdj$.  This process is summarized in
Algorithm \ref{ALGdiagonalSolve}.  Throughout the algorithms discussed in this
section, we treat the indices $s$ of blocks in $\bLdj$ as the labels of nodes
in an in-order traversal of a complete binary tree.  When we refer to a
{\em child} or {\em parent} of $s$, we mean the in-order index of the node
which is the child or parent of the node with in-order index $s$ in this tree.
For example, if $\bLdj$ has 7 blocks (as in \eqref{eq:ldjLabel}), then
$s = 4$ is the root of this tree and has left and right children with indices
3 and 6, respectively.

\begin{algorithm}[!t]

\SetKwData{Ii}{i}
\SetKwData{Dd}{d}
\SetKwData{Each}{each}
\SetKwData{Transpose}{transpose}
\SetKwInOut{Input}{input}
\SetKwInOut{Output}{output}

\Input{$j$, $\bG$, $\Transpose$, $s$}

\Begin{
  \SetAlgoVlined
  \CommentSty{// Leaf nodes correspond to dense diagonal blocks} \\[-2pt]
  \If{Node $s$ is a leaf} {
    \If{\Transpose} {
      \Return $(\dBlock{s})^{-T} \bG$
    } \Else {
      \Return $(\dBlock{s})^{-1} \bG$
    }
  } \Else {
    \CommentSty{// Partition $\bG$ into two parts} \\[-2pt]
    $\bG_{1} = \bG(1:|\cDs|,:)$ \\
    $\bG_{2} = \bG(|\cDs|+1:end,:)$ \\[5pt]

    \CommentSty{// Child blocks of $s$} \\[-2pt]
    $s_{1} = $ left in-order child of $s$ \\
    $s_{2} = $ right in-order child of $s$ \\[5pt]

    \If{\Transpose} {
      \CommentSty{// Backward substitution} \\[-2pt]
      $\bG_{1} \longleftarrow \diagonalSolve(j, \bG_{1}, \Transpose, s_{1})$ \\
      $\bG_{2} \longleftarrow \bG_{2} - \vBlock{s} \left( (\uBlock{s})^{T}
        \bG_{1} \right)$ \\
      $\bG_{2} \longleftarrow \diagonalSolve(j, \bG_{2}, \Transpose, s_{2})$
    } \Else {
      \CommentSty{// Forward substitution} \\[-2pt]
      $\bG_{2} \longleftarrow \diagonalSolve(j, \bG_{2}, \Transpose, s_{2})$ \\
      $\bG_{1} \longleftarrow \bG_{1} - \uBlock{s} \left( (\vBlock{s})^{T}
        \bG_{2} \right)$ \\
      $\bG_{1} \longleftarrow \diagonalSolve(j, \bG_{1}, \Transpose, s_{1})$
    }

    \Return $\left( \begin{array}{c} \bG_{1} \\ \bG_{2} \end{array} \right)$
  }
}

\cprotect\caption{{$\diagonalSolve$:
  Solves a system of equations using the submatrix of $\bLdj$ rooted at
  block $s$ in $\bLdj$'s hierarchical structure.
  }}
\label{ALGdiagonalSolve}
\end{algorithm}

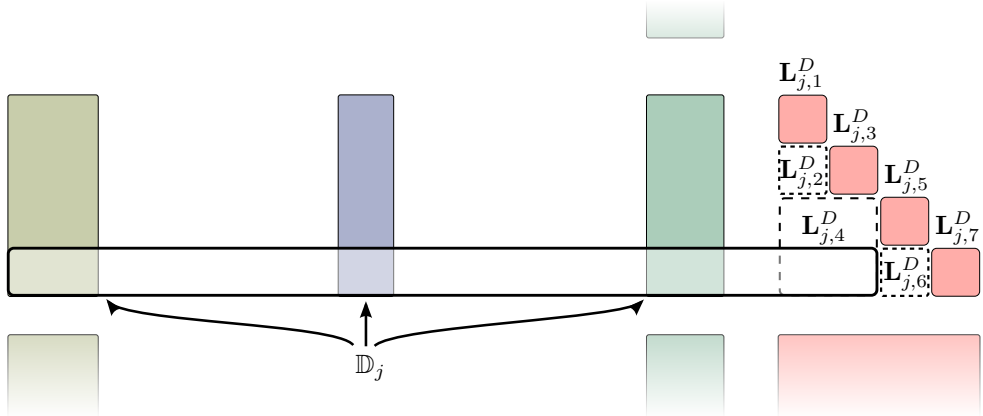
\begin{figure}
  \vspace{-6mm}
  \begin{tikzpicture}
  \node [anchor=south west] (label) at (0,0)
    {
      \resizebox{\hsize}{!}{%
      \includegraphics{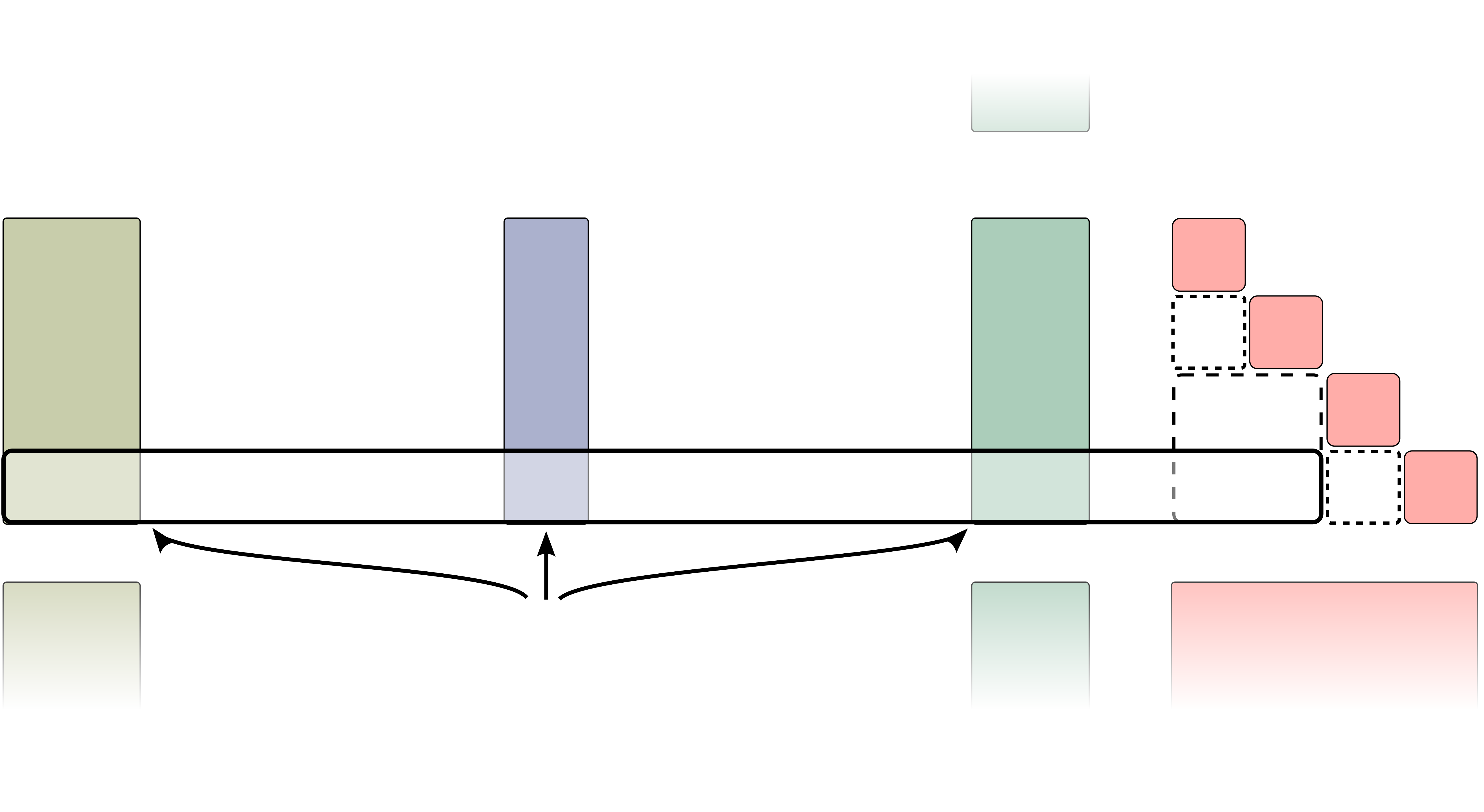}
      }%
    };
  \draw (4.975,1.7) node {$\bbD_{j}$};
  \draw (10.7,5.6) node {$\dBlock{1}$};
  \draw (11.42,4.92) node {$\dBlock{3}$};
  \draw (12.10,4.24) node {$\dBlock{5}$};
  \draw (12.78,3.56) node {$\dBlock{7}$};
  \draw (10.71,4.30) node {$\dBlock{2}$};
  \draw (12.09,2.95) node {$\dBlock{6}$};
  \draw (11.00,3.56) node {$\dBlock{4}$};
  \end{tikzpicture}
  \caption{Block row of a factor in which node $j$ (on the right) has a
           hierarchically compressed diagonal matrix
           $\bLdj$ with structure given by \eqref{eq:ldjLabel}.
           Compressed matrices are indicated with a dashed border.
           Consider the compressed block $\dBlock{6}$.
           When forming
           products with this matrix (for the purpose of compression), we
           require contributions from all columns of $\bL$ with non-zeros in
           the region highlighted with a thick black border.  Observe that
           this includes contributions from descendents $\bbD_{j}$, as well
           as other compressed blocks in $\bLdj$.  In this case, $\dBlock{6}$
           depends on three of node $j$'s descendents, as well as $\dBlock{4}$.
           }
  \label{FIGdiagDescendents}
\end{figure}

We will now use the $\diagonalSolve$ function to build an algorithm for
compressing an off-diagonal block $\dBlock{s}$.
As before, we multiply $\dBlock{s}$ with 
random matrices without explicitly constructing $\dBlock{s}$
As in \S\ref{SECodcompalg}, these operations 
depend on the contents of node $j$'s descendants $\bbD_{j}$.
In addition, we may need 
to consider 
contributions from previously
compressed blocks within $\bLdj$.
See Figure~\ref{FIGdiagDescendents} for a visual representation of this
dependence.

\begin{algorithm}[!h]

\SetKwData{Ii}{i}
\SetKwData{Dd}{d}
\SetKwData{Each}{each}
\SetKwData{Transpose}{transpose}
\SetKwData{True}{true}
\SetKwData{False}{false}
\SetKwData{Not}{not}
\SetKwData{Is}{is}
\SetKwData{Null}{null}
\SetKwInOut{Input}{input}
\SetKwInOut{Output}{output}

\Input{$j$, $\bG$, $\Transpose$, $s$}
\Output{$\vBlock{s}, \uBlock{s}$}

\Begin{
  \SetAlgoVlined
  \CommentSty{Apply necessary diagonal block inverse to $\bG$} \\[-2pt]
  $s_{1} = $ left child of $s$ \\
  $\bG \longleftarrow \diagonalSolve(j, \bG, \Transpose = \True, s_{1})$ \\[5pt]

  \CommentSty{// Initialize a workspace for multiplication} \\[-2pt]
  $\bW \longleftarrow \bA(\srDs, \scDs) \bG$ \\[5pt]

  \CommentSty{// Apply contributions from supernode descendants} \\[-2pt]
  \For{$\Each \, k \in \bbD_{j}$} {
    \CommentSty{// Extract the required sub-matrix from $\bG$} \\[-2pt]
    $\bG_{sub} \longleftarrow \gatherRows(\bG, \scDs, \srDkToj \cap \scDs)$
      \\[5pt]

    \CommentSty{// Form the needed product with two multiplications} \\[-2pt]
    $\bT \longleftarrow \left[ \bLok(\rDkToj \cap \cDs, :) \right]^{T}
      \bG_{sub}$ \\
    $\bT \longleftarrow \bLok(\rDkToj \cap \rDs, :) \bT$ \\[5pt]

    \CommentSty{// Scatter result to the output matrix} \\[-2pt]
    $\bW \longleftarrow \bW - \scatterRows(\bT, \srDkToj \cap \srDs, \srDs)$
  }
  $\,$ \\[-2pt]

  \CommentSty{// Apply contributions from other blocks in $\bLdj$} \\[-2pt]
  $p = $ parent of $s$ \\
  \While{$p $ \Is \Not \Null} {
    \CommentSty{// Block $p$ contributes to $s$ only if it appears}
      \\[-2pt]
    \CommentSty{// earlier in the ordering} \\[-2pt]
    \If{$p < s$} {
      \CommentSty{// Get necessary row and column ranges from $p$} \\[-2pt]
      \CommentSty{// This is valid since $\rDs \subset \rD{p}$ and
        $\cDs \subset \rD{p}$} \\[-2pt]
      $R_{sub} = \alignSet(\rDs, \rD{p})$ \\
      $C_{sub} = \alignSet(\cDs, \rD{p})$ \\[5pt]

      \CommentSty{// Perform multiplication similar to \eqref{eq:lowRankMult}}
        \\[-2pt]
      $\bU_{prod} \longleftarrow (\uBlock{p})^{T} \uBlock{p}$ \\
      $\bT \longleftarrow [\vBlock{p}(C_{sub},:)]^{T} \bG$ \\
      $\bT \longleftarrow \bU_{prod} \bT$ \\
      $\bT \longleftarrow \vBlock{p}(R_{sub},:) \bT$ \\[5pt]

      \CommentSty{Accumulate result in workspace} \\[-2pt]
      $\bW \longleftarrow \bW - \bT$ \\[5pt]
    }

    \CommentSty{Continue moving up the block hierarchy} \\[-2pt]
    $p \longleftarrow $ parent of $p$
  }

  \Return $\bW$
}

\cprotect\caption{{$\diagonalMultiply$:
  Given an input matrix $\bG$, forms the product $\dBlock{s} \bG$, where
  $\dBlock{s}$ is an off-diagonal block in $\bLdj$.  This is done without
  forming $\dBlock{s}$ explicitly.
  The $\alignSet$ function performs the following operation:
  $\alignSet(S_{1}, S_{2}) = \{ i - \min(S_{2}) + 1 : i \in S_{1} \}$.
  We will use this function to express row and column sets for block $s$
  relative to other blocks in the hierarchy.
  }}
\label{ALGdiagBlockMult}
\end{algorithm}

The $\diagonalMultiply$ algorithm (Algorithm \ref{ALGdiagBlockMult}) provides
the details of this procedure.  
Lines 7-15 
in this algorithm
resemble the descendant multiplication from Algorithm~\ref{ALGoffDiagMult}.
Lines 17-35 compute contributions from
other blocks inside of $\bLdj$.  As before, we also require the algorithm
$\diagonalMultiplyTrans$.  This algorithm has a similar structure to
Algorithm~\ref{ALGdiagBlockMult} 
With these two functions, 
we define a function
$\approximateDiagonalBlock$ which computes a low-rank representation of
$\dBlock{s}$ given some prescribed rank.  The structure of this algorithm
is not given here since it
is nearly identical to Algorithm~\ref{ALGapproxOffDiag}.

Finally, we turn to the question of how to construct dense diagonal blocks
within $\bLdj$.  As in Algorithm \ref{ALGfactorSupernode}, lines 7 \& 10, we
will consider update matrices built from off-diagonal blocks in node $j$'s
descendants.  In addition, we will need to consider contributions from
previously computed low-rank blocks in $\bLdj$.  The details of this process
are given in Algorithm \ref{ALGfactorDiagonal}.

\begin{algorithm}[!h]

\SetKwData{Ii}{i}
\SetKwData{Dd}{d}
\SetKwData{Each}{each}
\SetKwData{Transpose}{transpose}
\SetKwData{DiagUpdate}{diagUpdate}
\SetKwData{True}{true}
\SetKwData{False}{false}
\SetKwData{Not}{not}
\SetKwData{Is}{is}
\SetKwData{Null}{null}
\SetKwInOut{Input}{input}
\SetKwInOut{Output}{output}

\Input{$j$, $s$}
\Output{$\dBlock{s}$}

\Begin{
  \SetAlgoVlined
  \CommentSty{// Initialize Schur complement with matrix contents.} \\[-2pt]
  \CommentSty{// This is a diagonal block, so $\srDs = \scDs$.} \\[-2pt]
  $\schurBlock{s} \longleftarrow \bA(\srDs, \scDs)$ \\[5pt]

  \CommentSty{// Accumulate contributions from descendants.} \\[-2pt]
  \For{$\Each \, k \in \bbD_{j}$} {
    \CommentSty{// Build dense update block} \\[-2pt]
    $\DiagUpdate \longleftarrow \bLdk(\rDkTojs,:) \bLdk(\rDkTojs,:)^{T}$ \\[5pt]

    \CommentSty{// Scatter updates to the Schur complement} \\[-2pt]
    $\schurBlock{s} \longleftarrow \schurBlock{s}
      - \scatter(\DiagUpdate, \srDkTojs, \scDs, \srDkTojs, \scDs)$
  }
  $\,$ \\[-2pt]

  \CommentSty{// Accumulate contributions from previous blocks}
    \\[-2pt]
  \CommentSty{// in $\bLdj$.}
    \\[-2pt]
  $p \longleftarrow $ parent of $s$ \\
  \While{$p $ \Is \Not \Null} {
    \CommentSty{// Block $p$ contributes to $s$ only if it appears earlier}
      \\[-2pt]
    \CommentSty{// in the ordering} \\[-2pt]
    \If{$p < s$} {
      \CommentSty{// Get necessary row and column ranges from $p$} \\[-2pt]
      \CommentSty{// This is valid since $\rDs \subset \rD{p}$ and
        $\cDs \subset \rD{p}$} \\[-2pt]
      $R_{sub} = \alignSet(\rDs, \rD{p})$ \\
      $C_{sub} = \alignSet(\cDs, \rD{p})$ \\[5pt]

      \CommentSty{// Build a dense update matrix}
        \\[-2pt]
      $\bU_{prod} \longleftarrow (\uBlock{p})^{T} \uBlock{p}$ \\
      $\bT \longleftarrow \bU_{prod} [\vBlock{p}(C_{sub},:)]^{T}$ \\
      $\bT \longleftarrow \vBlock{p}(R_{sub},:) \bT$ \\[5pt]

      \CommentSty{// Subtract update from Schur complement} \\[-2pt]
      $\schurBlock{s} \longleftarrow \schurBlock{s} - \bT$
    }
  }

  \CommentSty{// Factor node $j$'s diagonal block} \\[-2pt]
  $\dBlock{s} \longleftarrow \Cholesky( \schurBlock{s} )$ \\[5pt]

  $\Return \, \dBlock{s}$
}

\cprotect\caption{{$\factorDiagonal$:
  Builds a factored diagonal block $\dBlock{s}$ within $\bLdj$.  $s$
  is assumed to be the index of a diagonal block.
  We introduce two new pieces of notation here:
  $\rDkTojs = \{ 1 \le p \le |\sR_{k}| : r_{k}^{p} \in \scDs \}$
  and
  $\srDkTojs = \{ r_{k}^{p} \in \sR_{k} : r_{k}^{p} \in \scDs \}$.
  That is, these are rows from descendant $k$ which are relevant to
  the formation of diagonal block $s$ within $\bLdj$. \\
  }}
\label{ALGfactorDiagonal}
\end{algorithm}

Given algorithms for forming diagonal and off-diagonal blocks in $\bLdj$,
building this matrix follows a straight forward process of iterating over the
blocks of $\bLdj$ in increasing order $s = 1, 2, \ldots$.  At each iteration,
we either form the diagonal block $\dBlock{s}$ or a low-rank decomposition
$\vBlock{s} (\uBlock{s})^{T}$.
Finally, we note that, unlike the off-diagonal compression scheme presented
in \S\ref{SECodblockcomp}, forming diagonal blocks and approximate off-diagonals
in this order does not guarantee that positive definiteness is maintained
throughout the factorization.  We discuss our simple method for addressing this
issue in \S\ref{SECavoidingIndefinite}.

\subsection{Choosing Diagonal Block Coordinates}
\label{SECdiagcoordinates}

In \S\ref{SECexposelowrank}-\ref{SECdiagcompalg} we assumed that indices
within $\bLdj$ could be reordered to expose low-rank structure in off-diagonal
blocks $\dBlock{s}$.  As we discussed in \S\ref{SECexposelowrank}, this is
accomplished by assigning spatial coordinates to degrees of freedom within
$\bLdj$.  Indices are reordered such that off-diagonal blocks in $\bLdj$
describe interactions between spatially separated ``pieces'' of the
separator with which node $j$ is associated.  However, we have not yet
discussed how these spatial coordinates are determined.  In many applications,
this information can be determined from the underlying PDE.  For example,
in \S\ref{SECresults} we discuss several model problems implemented in the
Deal.II finite element library.  For these problems, spatial coordinates
are determined directly from node positions in a finite element mesh.
Unfortunately, this information
may not be readily available in some cases.  In the interest of building a
general, algebraic preconditioner, we wish to also consider cases in which
spatial coordinates for system degrees of freedom are not provided.

Research in the area of {\em graph visualization} has yielded a variety of
methods for building visually appealing drawings of graphs
\cite{Kaufmann2001DGM376944,Battista98GraphDrawingAlgorithms}.  Given a sparse
matrix $\bA$, we can infer geometric positions for matrix indices by applying
these algorithms to the graph structure implied $\bA$'s non-zero pattern.
In this work, we appeal to {\em spectral graph drawing} algorithms, which
build positions based on the spectral properties of certain matrices associated
with the original system matrix $\bA$.  In particular, we consider $\bA$'s
{\em Graph Laplacian} $\bbL$, defined as follows:
\begin{equation}
  \bbL_{ij} =
    \left\{ \begin{array}{cl}
      -1 & \textrm{if } i \ne j, \, \bA_{ij} \ne 0 \\
      0 & \textrm{if } i \ne j, \, \bA_{ij} = 0 \\
      \left| \left\{ k \ne i : \bA_{ik} \ne 0 \right\} \right|
        & \textrm{if } i = j
    \end{array} \right.
\end{equation}
We evaluate the three lowest-order eigenvectors $\bv_{1}, \bv_{2}, \bv_{3}$
of $\bbL$ and associate the three-dimensional position
$[ \bv_{1i} \, \bv_{2i} \, \bv_{3i} ]$ with matrix index $i$.  We find that
low-accuracy approximations of these eigenvectors suffice, and we evaluate
these eigenvectors using Arnoldi iteration.  See \S\ref{SECresults} for
further discussion on the cost of constructing these positions.

\section{Additional Optimizations and Implementation Details}
\label{SECoptimizations}

In this section we discuss additional optimizations for further storage reduction
in our algorithm, as well as key implementation details. 

\subsection{Interior Blocks}
\label{SECinteriorBlocks}
The algorithm discussed in \S\ref{SECoffdiagcomp}-\ref{SECdiagcomp} builds
a sparse Cholesky factor on a matrix permuted with a nested dissection ordering.
Separators in the nested dissection
hierarchy that are sufficiently large -- that is, having more than $\tau_{O}$
variables -- are identified as supernodes in a
supernodal Cholesky factorization and the diagonal and off-diagonal blocks
for these supernodes are compressed.
Supernodes with fewer than $\tau_{O}$ variables may be factored using the
standard supernodal sparse Cholesky algorithm; however, in this section we
present a more memory-efficient method for handling these uncompressed blocks.

\begin{figure}
  \resizebox{\hsize}{!}{%
  \begin{tikzpicture}
  \node [anchor=south west] (label) at (0,0)
    {\includegraphics{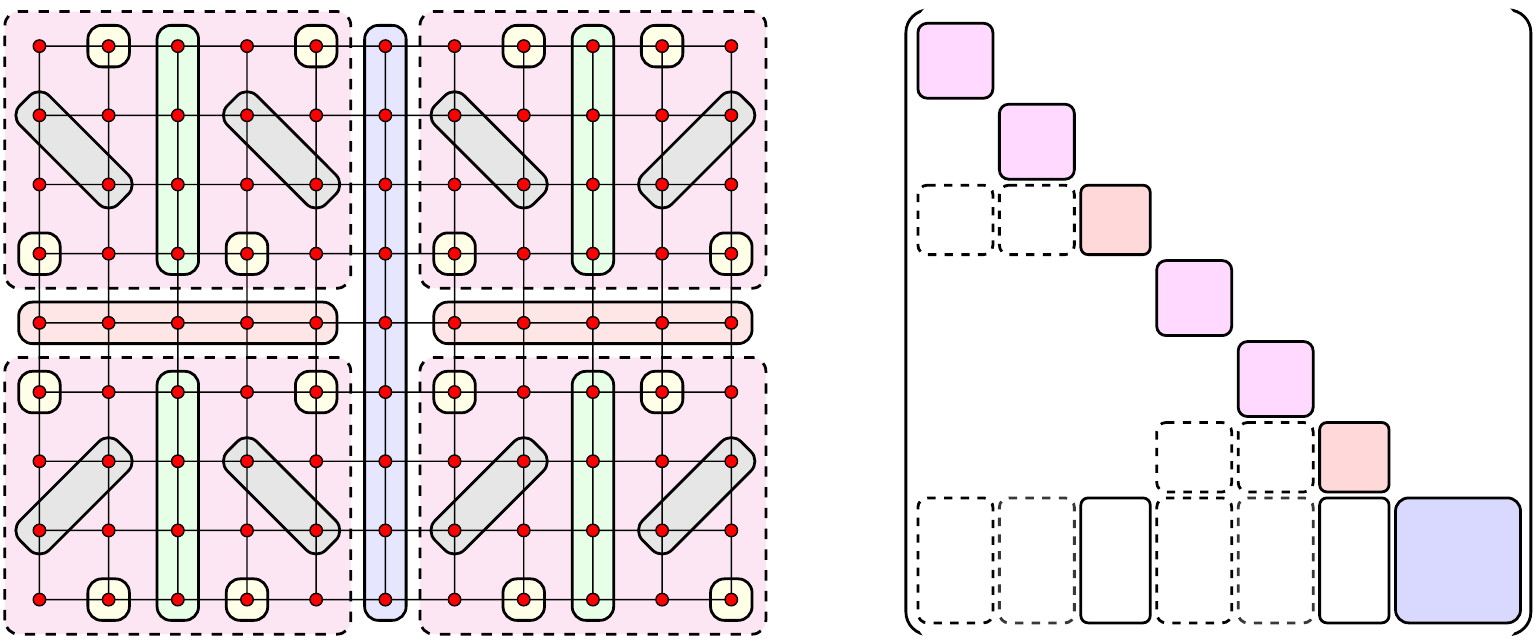}};
  \draw (-0.2,3.4) node {$j_{1}$};
  \draw (8.22,3.4) node {$j_{2}$};
  \draw (4.0, -0.1) node {$j_{3}$};
  \draw (-0.2, 0.5) node {$\mathbb{B}_{1}$};
  \draw (-0.2, 6.2) node {$\mathbb{B}_{2}$};
  \draw (8.22, 0.5) node {$\mathbb{B}_{3}$};
  \draw (8.22, 6.2) node {$\mathbb{B}_{4}$};

  \draw (11.43, -0.1) node {$j_{1}$};
  \draw (13.85, -0.1) node {$j_{2}$};
  \draw (14.90, -0.1) node {$j_{3}$};
  \draw (9.85, -0.1) node {$\mathbb{B}_{1}$};
  \draw (10.65, -0.1) node {$\mathbb{B}_{2}$};
  \draw (12.27, -0.1) node {$\mathbb{B}_{3}$};
  \draw (13.08, -0.1) node {$\mathbb{B}_{4}$};
  \end{tikzpicture}
  }%
  \caption{As in Figure~\ref{FIGoffdiagBlocks}, we apply nested dissection to
           a two-dimensional grid and set $\tau_{O} = 5$.  The separators
           with 5 or more variables partition this domain in to four disjoint
           subdomains, which we label as interior blocks $\bbB_{1,\ldots,4}$.
           These four blocks are shown as regions surrounded by a dashed line
           in the left image.
           The block structure in $\bL$ resulting from this partitioning is
           shown on the right.  As before, solid white blocks in this matrix
           denote off-diagonal blocks $\bLoj$ for which we apply compression
           according to the methods in \S\ref{SECoffdiagcomp}.  Dashed white
           blocks denote the off-diagonals of interior blocks
           $\bbB_{1,\ldots,4}$.  These blocks are not stored explicitly.
           We do, however, store the diagonal components of interior blocks
           explicitly, and these entries are evaluated using standard,
           supernodal Cholesky factorization restricted to the interior block
           subdomain.
          }
  \label{FIGinteriorBlocks}
\end{figure}
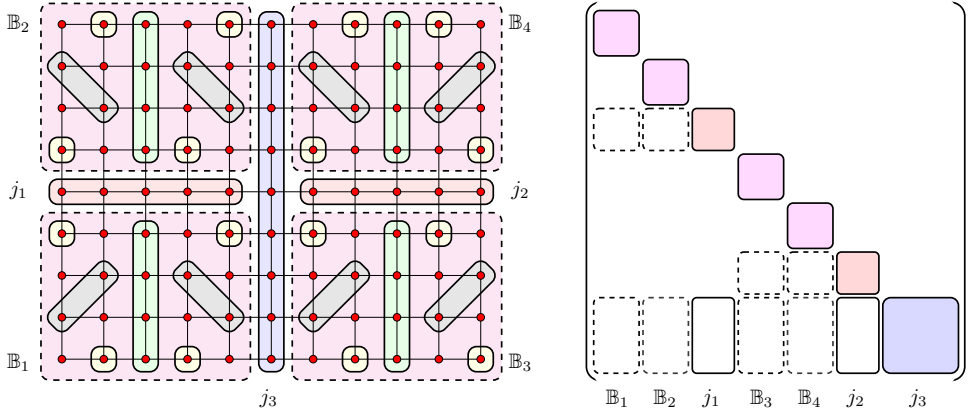

The collection of compressed
separators discussed above partitions the domain in to a collection of
mutually disjoint subdomains (see Figure~\ref{FIGinteriorBlocks} -- left side),
which we refer to as {\em interior blocks}.  This remains true even for
non-physical problems in which the ``domain'' is the graph defined
by the sparsity pattern of the matrix to be factored.
The nested dissection permutation guarantees that
the variables in an interior block appear in a contiguous block in the
reordered matrix (see Figure \ref{FIGinteriorBlocks} -- right side).
As such, an interior block can be represented by a sequential list of supernode
indices.  For interior block $i$ (numbered in the order in
which it appears in the reordered matrix), we use the notation
$\bbB_{i}$ to denote the list of supernode indices comprising the
block.
The column list $\sCbi$ and
off-diagonal row pattern $\sRbi$ associated with interior block $i$ are
defined as follows:
\begin{equation}
  \sCbi = \bigcup_{j \in \bbB_{i}} \sC_{j}
    \quad \quad \quad
  \sRbi = \left( \bigcup_{j \in \bbB_{i}} \sR_{j} \right)
    \setminus \sCbi
\end{equation}
The matrices $\bL(\sCbi, \sCbi)$ and $\bL(\sRbi, \sCbi)$ are the
diagonal and off-diagonal factor blocks for interior block $i$ (the dashed/white
and shaded matrix blocks in Figure \ref{FIGinteriorBlocks}, respectively).
Indices within interior block
$i$ are reordered via nested dissection to guarantee that
$\bL(\sCbi, \sCbi)$ is as sparse as possible.  However, for our purposes
we can think of interior blocks as representing the ``bottom'' level of
the nested dissection hierarchy.  When building a Cholesky factorization,
blocks in a nested dissection hierarchy only depend on blocks from lower levels
in this hierarchy.  Therefore, the factor contents for interior block $i$ are
evaluated as follows:
\begin{align}
  \bL(\sCbi, \sCbi) &= \chol(\bA(\sCbi, \sCbi))
    \label{EQinteriorDiagonalDef} \\
  \bL(\sRbi, \sCbi) &= \bA(\sRbi, \sCbi) \bL(\sCbi, \sCbi)^{-T}
    \label{EQinteriorOffDiagDef} 
\end{align}
We compute $\bL(\sCbi, \sCbi)$ using a standard, uncompressed supernodal
factorization.

The matrix $\bL(\sRbi, \sCbi)$ only stores interactions between the variables of
interior block $i$ and supernodes compressed using the methods of
\S\ref{SECoffdiagcomp}-\ref{SECdiagcomp} (see Figure~\ref{FIGinteriorBlocks}).
As a result, $\bL(\sRbi, \sCbi)$ may have many non-zero entries, making it
expensive to store explicitly.  Fortunately, we can 
still approximately factor $\bA$ without ever explicitly forming the block
$\bL(\sRbi, \sCbi)$.
In the standard supernodal factorization (Algorithm~\ref{ALGfactorSupernode}),
we explicitly form the Schur complement matrix for each supernode.
Here, $\bL(\sRbi, \sCbi)$ must be stored explicitly because it is required
when running Algorithm~\ref{ALGfactorSupernode} on ancestors of nodes in
interior block $i$.  The key insight of our approach is that our factorization
algorithm only uses $\bL(\sRbi, \sCbi)$ in places:
\begin{enumerate}
  \item We explicitly build Schur complements for {\em small} diagonal blocks
  with fewer than $\tau_{D}$ rows/columns.  This may depend on contributions
  from $\bL(\sRbi, \sCbi)$.

  \item We must be able to form products of the form
  $\bL(\sRbi, \sCbi) \bG$ and $\bL(\sRbi, \sCbi)^{T} \bG$ where $\bG$ is an
  arbitrary dense matrix.  These products are required by Algorithms
  \ref{ALGoffDiagMult} and \ref{ALGdiagBlockMult}.
\end{enumerate}

Forming the diagonal blocks referred to in the first requirement only
necessitates the formation of a sub-block of $\bL(\sRbi, \sCbi)$ with at
most $\tau_{D}$ rows (see Algorithm \ref{ALGfactorDiagonal}, line 8).  This
sub-block is needed exactly once for the formation of a diagonal block.
We can form small sub-blocks of $\bL(\sRbi, \sCbi)$ as needed to build diagonal
block Schur complements, then discard them immediately afterwards.

We also observe that the products from the second requirement listed above
can be formed without explicitly forming any part of $\bL(\sRbi, \sCbi)$.
Suppose that we wish to compress blocks in supernode $j$, and that this node
has some descendents in interior block $i$; that is,
$\bbD_{j} \cap \bbB_{i} \ne \emptyset$.  Forming the product $\bLoj \bG$ for
some dense matrix $\bG$ can be done by considering each descendent in $\bbD_{j}$
individually, as is done in Algorithm \ref{ALGoffDiagMult}.  Alternately, we
can consider all of the descendents in interior block $i$ simultaneously.
To do this, we first recall that $\bLbi = \bL(\sCbi, \sCbi)$ is computed
explicitly using sparse supernodal factorization, meaning that its inverse
can be applied quickly.  It follows from Algorithm \ref{ALGfactorSupernode}
that interior block $i$'s contribution to the Schur complement $\sUoj$ is
given by
\begin{equation}
\label{eq:interiorBlockSchur}
  \bL(\sR_{j}, \sCbi) \bL(\sC_{j}, \sCbi)^{T}.
\end{equation}
We can use (\ref{EQinteriorDiagonalDef}-\ref{EQinteriorOffDiagDef}) to
express the two matrices involved in \eqref{eq:interiorBlockSchur} as
\begin{equation}
  \bL(\sR_{j}, \sCbi) = \bA(\sR_{j}, \sCbi) \left( \bLbi \right)^{-T}
    \quad \quad \quad
  \bL(\sC_{j}, \sCbi) = \bA(\sC_{j}, \sCbi) \left( \bLbi \right)^{-T}
    \label{eq:interiorBlockClosedForm}
\end{equation}
Finally, we can use \eqref{eq:interiorBlockClosedForm} to write the product
of \eqref{eq:interiorBlockSchur} and an arbitrary dense matrix $\bG$
(as required by Algorithm \ref{ALGoffDiagMult}) as
\begin{equation}
  \bL(\sR_{j}, \sCbi) \bL(\sC_{j}, \sCbi)^{T} \bG
  =
  \bA(\sR_{j}, \sCbi) \left[ \left( \bLbi \right)^{-T}
  \left[ \left( \bLbi \right)^{-1} \left[ \bA(\sC_{j}, \sCbi)^{T} \bG \right]
  \right] \right]
  \label{eq:interiorBlockMult}
\end{equation}
As suggested by the parenthesis in \eqref{eq:interiorBlockMult}, this product
is the result of multiplying a sparse matrix with $\bG$, followed by two
sparse triangular solves involving $\bLbi$, followed by another multiplication
with a sparse matrix.  Since each of these operations can be carried out
efficiently, this provides an effective method for forming the matrix products
required by Algorithm~\ref{ALGoffDiagMult} without
having to explicitly store blocks of the form $\bL(\sRbi, \sCbi)$.
A similar method can be used to form products with blocks $\dBlock{s}$ by
replacing $\sR_{j}$ and $\sC_{j}$ in \eqref{eq:interiorBlockMult} (see
Algorithm \ref{ALGdiagBlockMult}) with $\srDs$ and $\scDs$, respectively.

While the optimizations discussed here have the potential to significantly
reduce storage requirements, this comes at the cost of somewhat more expensive
factorization and triangular solves.  We provide concrete examples of this
time-memory tradeoff in \S\ref{SECcomparisons}.

\subsection{Estimating Rank}
\label{SECestimatingRank}

\begin{figure}
\centering
\includegraphics[width=3.0in]{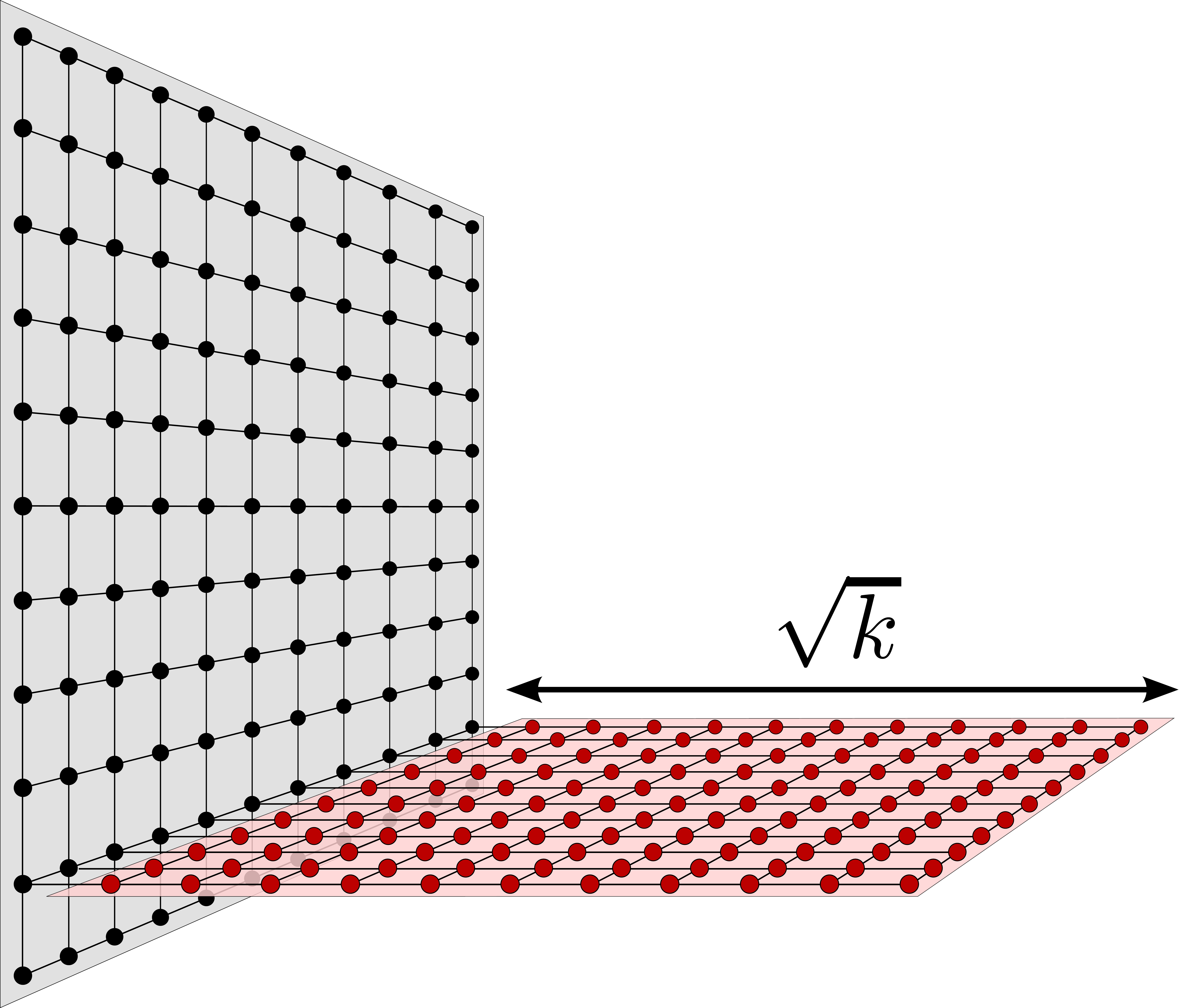}
\caption{Interaction between two square separators with grid-based topology.
         Assuming that each separator has $k$ variables, then the number of
         variables immediately adjacent to each other in the interaction
         between these two separators is $\sqrt{k}$.
         }
\label{fig:intersectingGrids}
\end{figure}

Up until now, we have assumed when building a low-rank approximations to
blocks in $\bL$ that the desired rank for each block is known {\em a priori}.
In this section, we discuss how block ranks are chosen.
We use our approximate rank-structured Cholesky factor as a preconditioner
for the
Preconditioned
Conjugate Gradient method.  Therefore, we are also free to use simple heuristics
to determine block ranks, with the understanding that the accuracy with which
we approximate blocks will influence the effectiveness of our preconditioner.
In principle, we could adaptively approximate blocks up to a certain
tolerance (see, e.g., \cite{liberty07,Halko2011}); however,
our experiments showed that the additional cost of
adaptive approximation outweighed the improved accuracy of
preconditioners built with this method.  Instead, we 
use a simple heuristic function depending only on the number of rows and
columns in the block to be approximated.  In particular, we assign the
following rank to a block $\bB \in \bbR^{m \times n}$:
\begin{equation}
  \blockRank(\bB) = \alpha \sqrt{k} \log_{2}( k ) + p
  \label{eq:rankFunction}
\end{equation}
where $k = \min(m, n)$ and $p$ is a small oversampling parameter.
To provide a brief, intuitive explanation as to why this function was
chosen, we consider
the interaction between two separators in 3D space.  Suppose that the separators
take the shape of regular, square, two-dimensional grids intersecting at
a right-angle (see Figure \ref{fig:intersectingGrids}).  This may be the
case in, for instance, a PDE discretized on a regular, three-dimensional finite
difference grid.  Assuming that each separator has $k$ variables, there are
$\sqrt{k}$ in each separator which are immediately adjacent to the other
separator.  It is for this reason that we include a term proportional to
$\sqrt{k}$ in \eqref{eq:rankFunction}.  We also scale $\blockRank(\bB)$
by $\log_{2}(k)$ since we empirically observe better preconditioning behavior
when larger ranks are used to approximate larger blocks from $\bL$.  In
practice, we also use two different constants in \eqref{eq:rankFunction}
-- $\alpha^{D}$ and $\alpha^{O}$ -- which determine ranks during diagonal and
off-diagonal compression, respectively.

\subsection{Avoiding Indefinite Factorizations}
\label{SECavoidingIndefinite}

As discussed originally in \S\ref{SECdiagcompalg}, our scheme for compressing
diagonal blocks $\bLdj$ does not provide a guarantee that Schur complements
formed during the factorization will remain positive definite.  This could
result in our factorization algorithm failing for certain inputs.
Fortunately, the argument in \S\ref{SECodblockcomp} guarantees that in the
absence of compression of diagonal blocks $\bLdj$, all Schur complements
remain positive definite.  This implies that we can avoid the indefinite
Schur complements by approximating diagonal blocks with sufficient
accuracy.  In practice, we address this issue by adapting the diagonal
compression parameter $\alpha_{D}$ in the event that an indefinite diagonal
matrix is encountered during factorization.  Specifically, we initialize
$\alpha_{D} = 0.5$ and if factorization fails due to an indefinite matrix,
we increase this constant $\alpha_{D} \longleftarrow 1.25 \alpha_{D}$ and
restart the factorization process.  This strategry increases the accuracy
with which diagonal blocks are approximated until factorization is successful.

\section{Results}
\label{SECresults}

We have applied the method described in this paper to a number of sample
problems.
In \S\ref{sec:nonlineardescription} we demonstrate the behavior of our solver
on a challenging nonlinear elasticity problem.  When applied to the linear
systems arising in this problem, our solver provides significant performance
improvements over a variety of standard solvers.  In \S\ref{sec:otherproblems},
we discuss the behavior of our solver on a variety of other sample problems.
We consider both standard examples implemented using the
{\em deal.II} finite element analysis library
\cite{BangerthHartmannKanschat2007,DealIIReference} and examples taken from
the {\em University of Florida sparse matrix collection} \cite{Davis:2011}.
While the performance differences between our solver and standard
direct and iterative solvers are less dramatic in these examples, these
results demonstrate the robustness of our solver.

\subsection{An Example: A Nonlinear Elasticity Problem}
\label{sec:nonlineardescription}

In this section, we discuss an example that illustrates the behavior
of our current solver.  While we have tested our solver on numerous
problems, the problem described here poses particular difficulty for
standard iterative methods.  As such, it is an ideal candidate for a
hybrid approach such as ours which leverages the reliability of direct
solvers with the low memory overhead of iterative methods.

We evaluate our problem using a benchmark problem taken
from~\cite{reese00} based on a standard example from the deal.II
finite element analysis
library~\cite{DealIIReference,BangerthHartmannKanschat2007}.  This
simulation models quasi-static loading of a nearly-incompressible,
hyperelastic block under compression.  The code
\begin{wrapfigure}[18]{r}{2.75in}
\centering
\vspace{-3mm}
\includegraphics[width=2.75in]{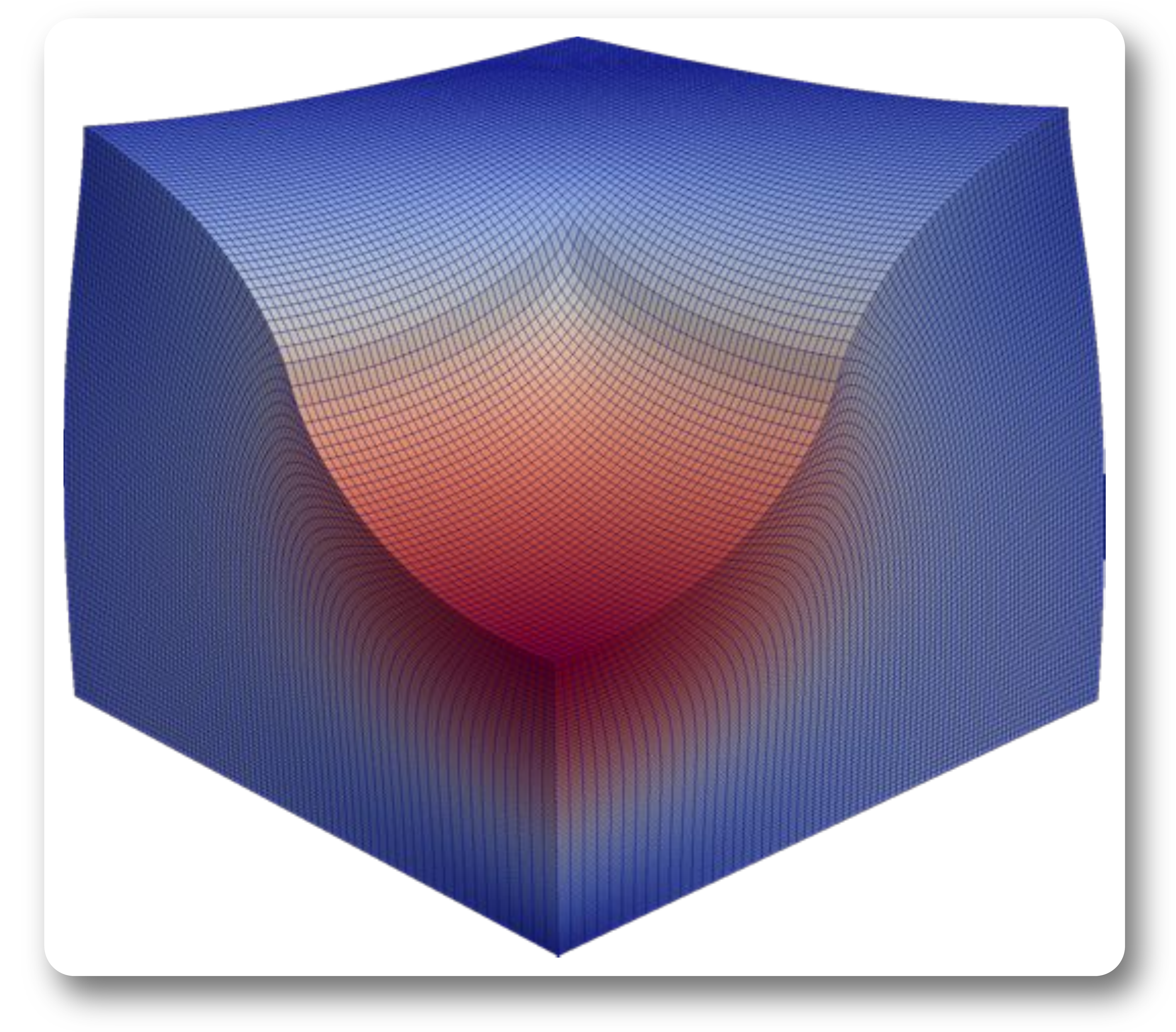}
\figspace{}
\figspace{}
\caption{High-resolution ($N = 80$), elastic block under compression.
         Our sparse rank-structured preconditioner was used to compute the
         deformations seen here.}
\label{fig:highResBlock}
\end{wrapfigure}
applies a force incrementally over two load steps;
at each load step, a nonlinear system of equations is solved to
determine the resulting deformation of the block.  The nonlinear
system is solved by a Newton iteration, and we evaluate the
performance of our solver for solving the sequence of linear systems
that arise during this process.



Because this problem is nearly incompressible, standard
displacement-based elements would be prone to locking.  Consequently,
our test problem uses a mixed formulation with explicit pressure and
dilation field variables in addition to the displacement fields.
We consider two versions of this problem -- henceforth referred to as the
$p=1$ and $p=2$ problems.
In the $p=1$ problem, displacements are discretized with continuous linear
Lagrange brick elements, while pressure and dilation are discretized using
discontinuous piecewise constant functions.
The $p=2$ problem discretizes displacements with quadratic elements and
uses discontinuous linear elements for pressure and dilation.
All variables are discretized on an $N \times N \times N$ element grid.
In all problem instances, the pressure and
dilation variables are condensed out prior to the linear solve, so the
system we solve involves only displacement variables.

\subsubsection{Comparison to Standard Iterative Solvers}
\label{sec:iter_perf}

\begin{figure}[ht]
  \begin{center}
    \subfigure
    {
      \label{fig:convergencePlotp150}
      \includegraphics{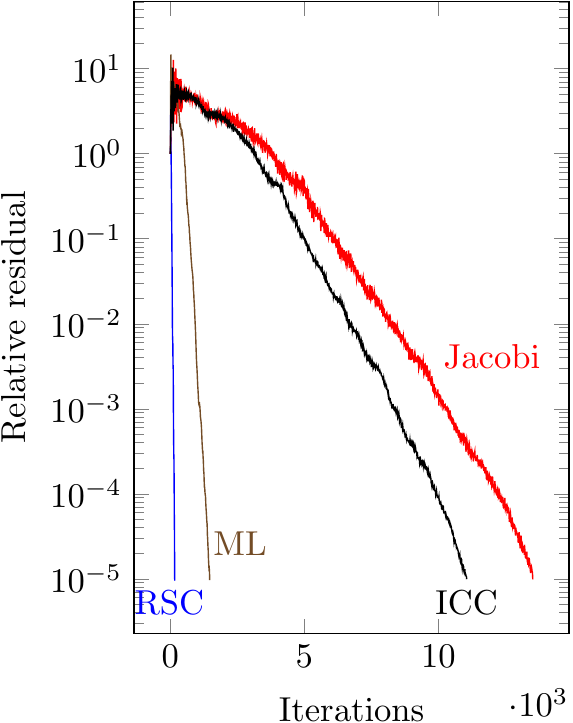}
    }
    \subfigure
    {
      \label{fig:accuracyVsTimep150}
      \includegraphics{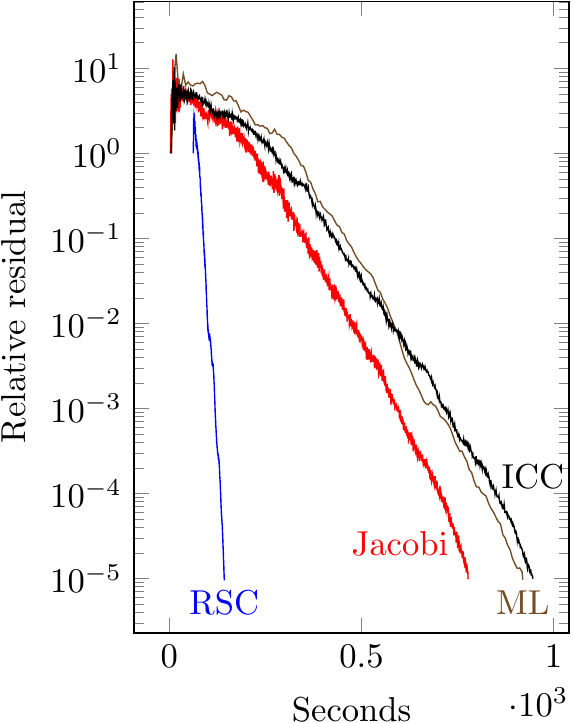}
    }
  \end{center}
\caption{PCG convergence for the first linear solve in the $p=1$ benchmark
  problem with $N = 50$.  Jacobi, ICC, ML and RSC refer to solves
  preconditioned with Jacobi, incomplete Cholesky (IFPACK), multigrid and
  rank-structured Cholesky (our solver) preconditioners, respectively.
  Convergence curves relative to wall clock time start at $t > 0$ due
  to time required to construct the preconditioner.}
\end{figure}

\begin{figure}[ht]
  \begin{center}
    \subfigure
    {
      \label{fig:convergencePlotp235}
      \includegraphics{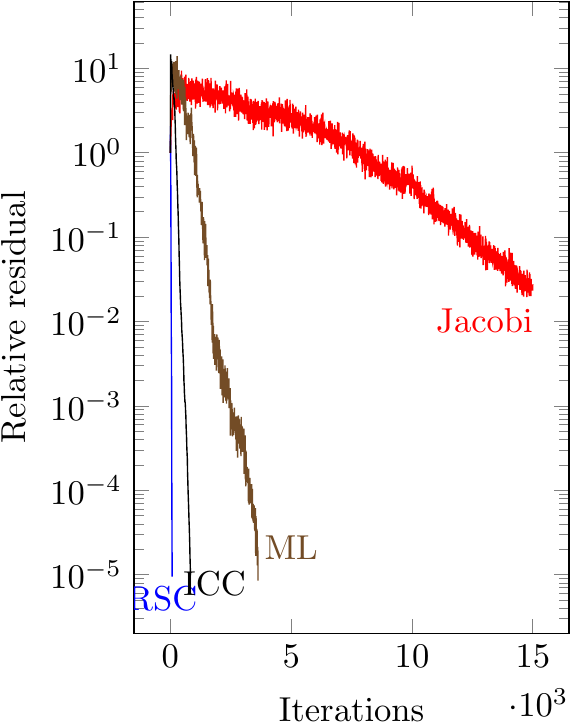}
    }
    \subfigure
    {
      \label{fig:accuracyVsTimep235}
      \includegraphics{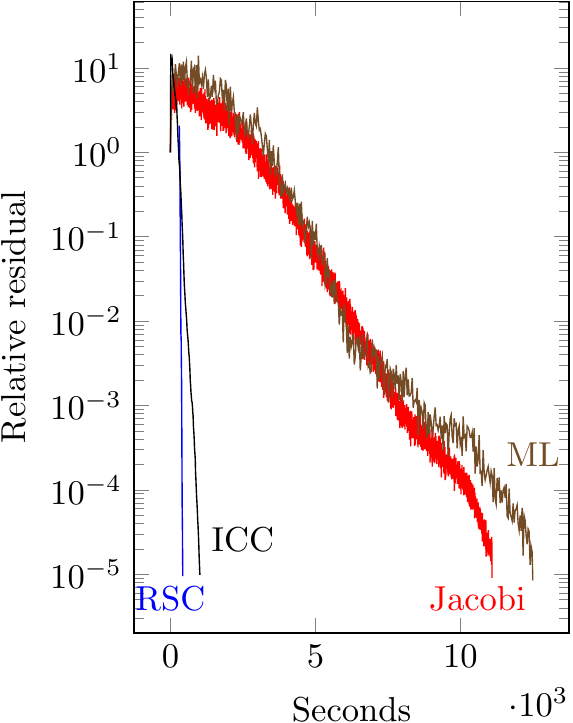}
    }
  \end{center}
\caption{PCG convergence for the first linear solve in the $p=2$ benchmark
         problem with $N = 35$ (see the caption for Figure
         \ref{fig:convergencePlotp150}-\ref{fig:accuracyVsTimep150}).
         Note that the Jacobi solve was run to convergence: the left
         plot is truncated for clarity.}
\end{figure}

\begin{table}
\begin{center}
   \resizebox{1.0\hsize}{!}{
   \begin{tabular}{c|r|c|c|c|c|c|c|c}
     \whline{0.8pt}
     \rule{0pt}{3mm}
     &$N$ & 20 & 30 & 40 & 50 & 60 & 70 & 80 \\
     &$n$                      & 27783    & 89373    & 206763   &
                                397953   & 680943   & 1073733  & 1594323\\
     \whline{0.6pt}
     \multirow{4}{*}{Jacobi}
     &Total time (s)     & 194      & 1170     & 4065      &
                                10520    &&&\\
     &Mean time (s)   & 14       & 84       & 271
                              & 701      &&&\\
     &Total iterations
        & 49579   & 87957   & 133361  & 170750 &&&\\
     &Mean iterations
        & 3541    & 6282    & 8890    & 11383  &&&\\
     \whline{0.6pt}
     &
     Total time (s)     & 207      & 1334     & 5159      &
                                13540    &&&\\
     ICC & Mean time (s)   & 15       & 95       & 344
                              & 903      &&&\\
     (IFPACK) & Total iterations
        & 48637   & 87244   & 132981  & 168172 &&&\\
     &Mean iterations
        & 3474    & 6231    & 8865    & 12110  &&&\\
     \whline{0.6pt}
     \multirow{4}{*}{ML}
     &Total time (s)     & 270      & 1447     & 4763     &
                                11360    &&&\\
     &Mean time (s)   & 19       & 103      & 318
                              & 757      &&&\\
     &Total iterations
        & 5949    & 10170   & 15249   & 19109  &&&\\
     &Mean iterations
        & 425     & 726     & 1017    & 1274   &&&\\
     \whline{0.6pt}
     \multirow{4}{*}{RSC} &
     Total time (s)     & 32       & 173      & 618      &
                                1470     & 2862     & 6174     & 12350
                                \\
     &Mean time (s)   & 2.3       & 12        & 41
                              & 98     & 191       & 412      & 823
                              \\
     &Total iterations
        & 349     & 615     & 1209    & 1514   & 1729   & 2278    & 4050
                  \\
     &Mean iterations
        & 24      & 41      & 80      & 100    & 115    & 151     & 270
        \\
     \whline{0.8pt}
   \end{tabular}
   }
\end{center}
\caption{{\bf Nonlinear elasticity ($p=1$) performance results:}
  Relative performance of standard iterative methods (Jacobi,
  incomplete Cholesky, and multigrid) compared to rank-structured
  Cholesky.
  }
\label{tbl:nonlinearComparison_p1}
\end{table}

\begin{table}
\begin{center}
   \resizebox{1.0\hsize}{!}{
   \begin{tabular}{c|r|c|c|c|c|c|c|c|c}
     \whline{0.8pt}
     \rule{0pt}{3mm}
     &$N$ & 10 & 15 & 20 & 25 & 30 & 35 & 40 & 45 \\
     &$n$                      & 27783    & 89373    & 206763   &
                                397953   & 680943   & 1073733  & 1594323 &
                                 \\
     \whline{0.6pt}
     \multirow{4}{*}{Jacobi}
     &Total time (s)     & 501      & 3731     & 13610      &
                                33920    &&&&\\
     &Mean time (s)   & 36       & 267       & 907
                              & 2261      &&&&\\
     &Total iterations
        & 60289   & 110971   & 175525  & 231833 &&&&\\
     &Mean iterations
        & 4306    & 7926    & 11701    & 15455  &&&&\\
     \whline{0.6pt}
     &
     Total time (s)     & 61      & 287     & 955      &
                                2473    & 5097    & 10150  &&\\
     ICC & Mean time (s)   & 4.4       & 20       & 64
                              & 165     & 340    & 677   &&\\
     (Aztec) & Total iterations
        & 897   & 1590   & 2689  & 3914 & 5341  & 7030  &&\\
     &Mean iterations
        & 64    & 113    & 179    & 260  & 356  & 468   &&\\
     \whline{0.6pt}
     \multirow{4}{*}{ML}
     &Total time (s)     & 918      & 4431     & 13920     &
                                34950    & 73690   & 144900 & 227800 &\\
     &Mean time (s)   & 66       & 317      & 928
                              & 2330      & 4913   & 9660   & 15187  &\\
     &Total iterations
        & 9149    & 13770   & 20748   & 26534  & 32292  & 40126  & 43509 &\\
     &Mean iterations
        & 653     & 983     & 1383    & 1768   & 2152   & 2675   & 2900  &\\
     \whline{0.6pt}
     \multirow{4}{*}{RSC} &
     Total time (s)     & 42       & 200      & 616      &
                                1460     & 2844     & 5139     & 9873   &
                                14900    \\
     &Mean time (s)   & 3       & 13        & 41
                              & 97      & 190       & 668      & 658    &
                              993       \\
     &Total iterations
        & 121     & 196     & 343     & 438    & 497    & 668     & 1438  &
        1429      \\
     &Mean iterations
        & 8       & 13      & 22      & 29     & 33     & 44      & 95    &
        95        \\
     \whline{0.8pt}
   \end{tabular}
   }
\end{center}
\caption{{\bf Nonlinear elasticity ($p=2$) performance results:}
  Relative performance of standard iterative methods (Jacobi,
  incomplete Cholesky, and multigrid) compared to rank-structured
  Cholesky.
  }
\label{tbl:nonlinearComparison_p2}
\end{table}

We have solved the benchmark problem with the preconditioned conjugate
gradient (PCG) iteration using our rank-structured Cholesky
preconditioner and several other preconditioners provided in
Trilinos~\cite{Trilinos-Overview,Trilinos-Dev-Guide,Trilinos-Dev-Guide-II,
  Trilinos-Users-Guide,1089021}, a library of high-performance solvers
developed primarily at Sandia national labs.  Our code consistently
out-performed a Jacobi preconditioner, an incomplete Cholesky (ICC)
preconditioner, and a multi-level (ML) preconditioner in both
iteration counts and wall clock time
(Figure~\ref{fig:convergencePlotp150}--\ref{fig:accuracyVsTimep150}).  Our
timing results are summarized in Table~\ref{tbl:nonlinearComparison_p1}.

All results reported in this section were generated on an 8-core Intel
Xeon X5570 workstation with 48GB of memory running Ubuntu 12.04, with
LAPACK and BLAS implementations provided by the Intel Math Kernel
Library version 11.0.  We use the preconditioned conjugate gradient
(PCG) implementation provided by AztecOO for all tests.  All linear
systems were solved to a relative $\ell_{2}$ residual error threshold
of $10^{-5}$.  At this accuracy level, the nonlinear iteration
required 14--15 linear solve steps to converge (as compared to 12
linear solves for a standard Cholesky solver).  Due to the time
required to solve linear systems using the standard preconditioners,
we only ran these example for $N \le 50$.  For the rank-structured
Cholesky solver, we ran examples up to $N \le 80$.


We observe the following properties for the different preconditioners
for this problem:

\myparagraph{Jacobi:} The Jacobi preconditioner (diagonal
preconditioner) is simple, but it usually only modestly accelerates
convergence.
When solving the $p=1$ problem, more {\em iterations} are required to converge
with the Jacobi preconditioner than with more sophisticated
preconditioners.  However, each iteration is so cheap that this process
requires less {\em wall clock time} than solves performed with other standard
preconditioners.  Meanwhile, in the $p=2$ problem, the number of iterations
required when using a Jacobi preconditioner grows considerably, making other
standard preconditioners more competitive.

\myparagraph{ICC:} We timed the PCG iteration using incomplete
Cholesky (ICC) preconditioners implemented in both IFPACK and AztecOO.
As with most
incomplete factorization codes, these solvers
require several parameters, including the level of fill allowed in
the factorization, drop tolerances dictating which matrix entries
should be discarded, and parameters controlling perturbations to the
matrix's diagonal.  The latter are required to avoid poorly
conditioned factorizations since solvers based on incomplete
factorizations appear to encounter severe conditioning issues when
applied to this problem.  In general, ``good'' parameter choices
depend on the problem.  We chose parameters for this problem based on
experiments with a small problem instance (e.g., $N = 20$ for the $p=1$
problem).  When solving the $p=1$ problem, even with the
tuned parameters, we require almost as many iterations with IFPACK's ICC
preconditioner as
with the much simpler Jacobi preconditioner.  Moreover, the incomplete
Cholesky preconditioner costs more than applying the Jacobi
preconditioner, so the overall time to solve the linear systems is
actually {\em larger} in this case.  Meanwhile, we find that Aztec's
ICC preconditioner is unable to make any progress towards convergence in
the $p=1$ problem.  We observe the opposite behavior in the
$p=2$ problem.  Specifically, IFPACK's ICC preconditioner makes no progress
towards convergence, whereas the AztecOO preconditioner performs reasonably well
(but still significantly slower than our rank structured solver).
In fact, while IFPACK's ICC solver was the least effective standard solver
(in terms of wall clock time) applied to the $p=1$ problem, AztecOO's ICC
solver was the {\em most} effective standard solver applied to the $p=2$
problem.  This phenomenon provides further evidence of the difficulty
associated with choosing an effective preconditioner for a given problem.

\myparagraph{ML:} 
While multigrid preconditioners perform well on many problems, on our
benchmark we see relatively poor convergence and long solve times.  We
use an algebraic multigrid preconditioner that solves the problem at
its coarsest level using a direct solver provided by Amesos.  The next
coarsest level applies a symmetric Gauss-Seidel smoother over several
sweeps (4 in this case).  Finer levels use a degree-2 Chebyshev
polynomial smoother.  These parameters were chosen based on good
convergence behavior on a small problem instance (e.g., $N = 20$ for the $p=1$
problem).
When solving the $p=1$ problem, this
solver requires significantly fewer iterations than the Jacobi or
incomplete Cholesky solvers; but because applying the preconditioner
is relatively expensive, it takes about as long to solve with the
multigrid preconditioner as with a Jacobi preconditioner.
When applied to the $p=2$ problem, the multigrid preconditioner requires
significantly fewer iterations than the Jacobi solver, but many more than
the Aztec00 incomplete Cholesky solver.  As a result, this solver is the
least effective standard solver that we tested on the $p=2$ problem.
We also note here that multigrid frameworks such as the one provided by ML
require
tuning a wide variety of parameters (some of which are discussed
above).  In some cases, tuning problem-specific parameters to achieve
good convergence behavior may outweigh the cost of solving the problem
with a simpler method.

\myparagraph{Rank-Structured Cholesky:} 
The conjugate gradient method preconditioned with our rank-structured
Cholesky solver converges quickly, both in terms of the iteration count
and in terms of wall clock time.  This is a significant improvement
over the other preconditioners.

\subsubsection{Comparison to Exact Factorization}
\label{sec:standard_perf}

\begin{figure}[ht]
  \begin{center}
    \subfigure{
      \label{fig:timingPlot}
      \begin{tikzpicture}
\begin{axis}[
  xlabel={$n$},ylabel={Solve time (s)},
  width=0.45\textwidth,
  height=8cm
]
\addplot coordinates {
  (27783,   3)
  (89373,   15)
  (206763,  47)
  (397953,  121)
  (680943,  242)
  (1073733, 549)
  (1594323, 1062)
} node[left] at (axis cs:1594323,1062) {RSC};
\addplot coordinates {
  (27783,   14)
  (89373,   84)
  (206763,  271)
  (397953,  701)
} node[right] at (axis cs:397953,701) {Jacobi};
\addplot coordinates {
  (27783,   19)
  (89373,   103)
  (206763,  318)
  (397953,  757)
} node[right] at (axis cs:397953,802) {ML};
\addplot coordinates {
  (27783,   15)
  (89373,   95)
  (206763,  344)
  (397953,  903)
} node[right] at (axis cs:397953,903) {ICC};
\addplot coordinates {
  (27783,   1.5)
  (89373,   7.7)
  (206763,  27)
  (397953,  76)
  (680943,  186)
  (1073733, 409)
} node[right] at (axis cs:1073733,409) {Cholesky};
\end{axis}
\end{tikzpicture}
    }
    \subfigure{
      \label{fig:memoryPlot}
      \begin{tikzpicture}
\begin{axis}[
  xlabel={$n$},ylabel={Memory (GB)},
  width=0.45\textwidth,
  height=8cm
]
\addplot coordinates {
  (27783,   0.053)
  (89373,   0.203)
  (206763,  0.475)
  (397953,  0.958)
  (680943,  1.830)
  (1073733, 3.074)
  (1594323, 4.700)
} node[above left] at (axis cs:1.6e6,4.7) {RSC};
\addplot coordinates {
  (27783,   0.237)
  (89373,   0.965)
  (206763,  3.187)
  (397953,  7.625)
  (680943,  16.044)
  (1073733, 30.003)
} node[right] at (axis cs:1.1e6,30) {Cholesky};
\end{axis}
\end{tikzpicture}
    }
  \end{center}
  \caption{Mean times and memory usage for the $p=1$ benchmark problem with
    problem between $N=20$ and
    $N=80$.  Solvers based on direct factorization are much faster
    than the standard preconditioners in this benchmark.  Although it
    takes somewhat longer to solve systems with our rank-structured
    solver than with an exact factorization, the memory requirements
    of the latter approach make it infeasible for larger problems.}
  \label{fig:memoryTimePlots}
\end{figure}
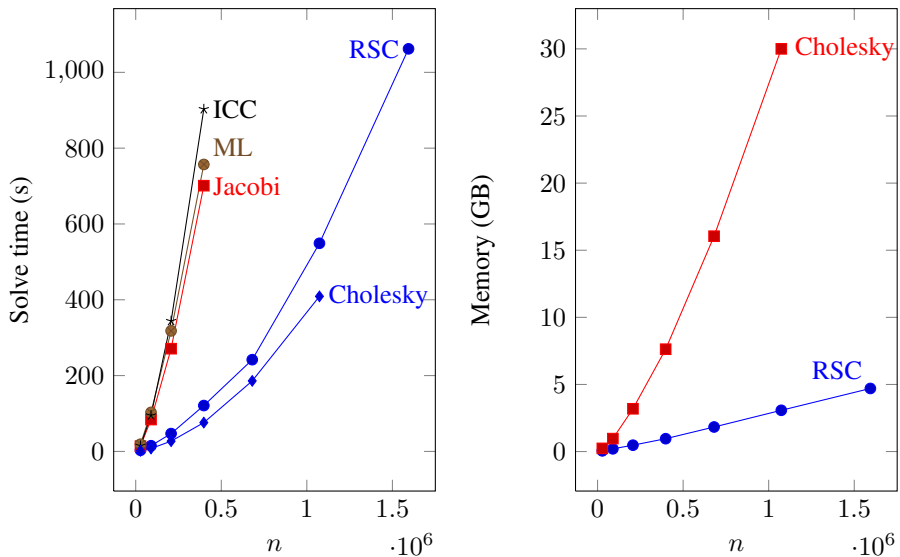

Beside comparing to standard preconditioners, we also compare our code
to an exact sparse Cholesky factorization.  As discussed above, our
sparse Cholesky implementation closely mirrors CHOLMOD~\cite{chen08}
and achieves similar performance and memory usage for exact
factorizations.  Because CHOLMOD places restrictions on problem size,
we use our code for both the rank-structured approximation and the
exact Cholesky factorizations

In Figure~\ref{fig:memoryTimePlots}, we show how much time and memory
we need to solve linear systems with the rank-structured and exact
sparse Cholesky factorizations.  For this benchmark, both the
rank-structured and the exact Cholesky solvers are much faster than
the standard preconditioned iterations.  The rank-structured Cholesky
solver is somewhat slower than the exact Cholesky solver; but the
memory requirements of the latter approach make it infeasible for
larger problems.

\subsection{Other Sample Problems}
\label{sec:otherproblems}
In this section, we discuss other sample problems arising either from
finite element discretizations in Deal.II or from the University of Florida
Sparse Matrix collection.  For some examples, we also consider the effect
of varying the number of power iterations used for low-rank approximation
in our solver (see \S\ref{SECodcompalg}).  We find that increasing this
number can result in somewhat more accurate approximation, and a modest
reduction in PCG iterations. \\

\myparagraph{Finite element analysis of a trabecular bone:}
We consider the stiffness matrix produced by finite element analysis of a
three-dimensional trabecular bone model.  The matrix used in this problem is
provided in the University of Florida sparse matrix collection
\footnote{\url{http://www.cise.ufl.edu/research/sparse/matrices/Oberwolfach/bone010.html}},
has dimension $n = 986703$ and has 24419243 non-zeros in it's lower
triangular component.  We use low-accuracy eigenvectors of the matrix's
graph Laplacian to compute three-dimensional coordinates for degrees of freedom
in this system (see \S\ref{SECdiagcoordinates}).  Since no right-hand-side
vector is provided for this problem, we solve the system $\bA \bx = \bb$ with
the constant vector $\bb = ( 1 \,\, 1 \,\, \ldots \,\, 1 )^{T}$.  \\

\begin{figure}[ht]
  \begin{center}
    \subfigure
    {
      \label{fig:convergencePlotBone}
      \includegraphics{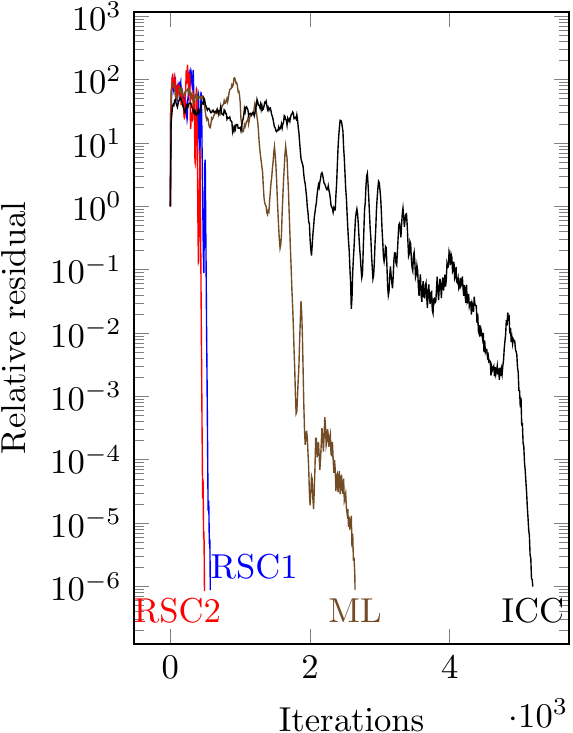}
    }
    \subfigure
    {
      \label{fig:accuracyVsTimeBone}
      \includegraphics{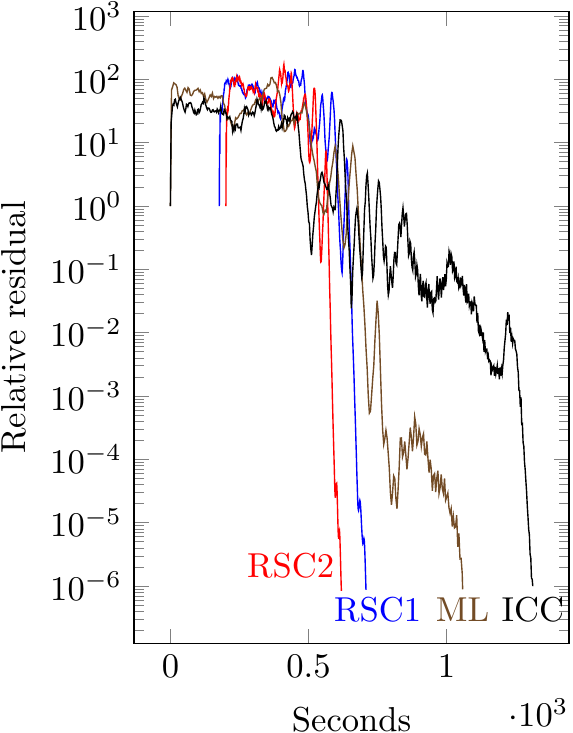}
    }
  \end{center}
\caption{PCG convergence for the trabecular bone problem.  We compare results
         using incomplete Cholesky (ICC), a multigrid preconditioner (ML), and
         our solver using $s=1$ or $2$ power iterations when building
         low-rank approximations (RSC1 and RSC2, respectively).  We see that
         the improved accuracy from more power iterations results in a modest
         reduction in solution time.}
\end{figure}

\myparagraph{Finite element analysis of a steel flange:}
Here, we consider a linear system arising from a three-dimensional
mechanical problem discretizing a steel flange.  This example can be found in
the University of Florida sparse matrix collection
\footnote{\url{http://www.cise.ufl.edu/research/sparse/matrices/Janna/Flan_1565.html}},
has dimension $n = 1564794$ and has 57865083 non-zeros in it's lower triangular
component. \\

\begin{figure}[ht]
  \begin{center}
    \subfigure
    {
      \label{fig:convergencePlotFlange}
      \includegraphics{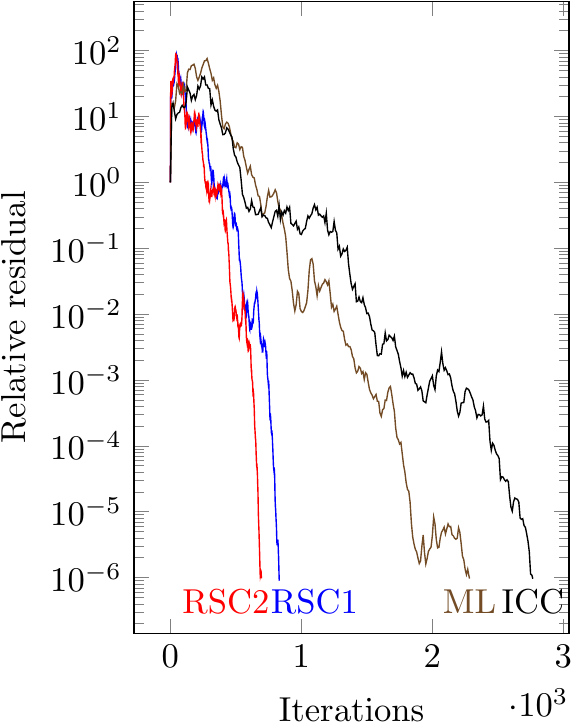}
    }
    \subfigure
    {
      \label{fig:accuracyVsTimeFlange}
      \includegraphics{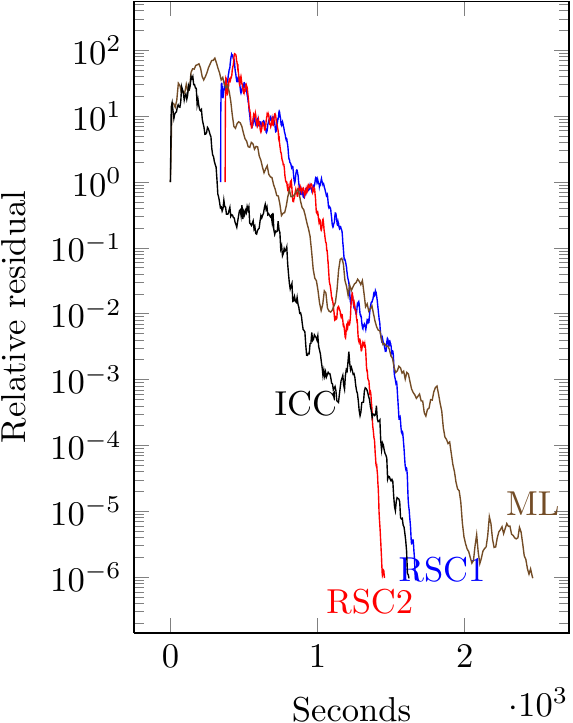}
    }
  \end{center}
\caption{PCG convergence for the steel flange problem.}
\end{figure}

\myparagraph{Poisson's equation:}
We apply our solver to Poisson's equation
\begin{align}
  - \nabla \cdot \bK(\bx) \nabla p &= f \textrm{ in } \Omega 
    \label{eq:poissonEq} \\
  p &= g \textrm{ on } \partial \Omega
    \label{eq:poissonBC}
\end{align}
where the coefficients $\bK(\bx)$ are optionally both inhomogeneous and
anisotropic.  The basic setup of this problem follows a standard example
from the deal.II library
\footnote{\url{http://www.dealii.org/developer/doxygen/deal.II/step_20.html}}.
In particular, rather than solving (\ref{eq:poissonEq}-\ref{eq:poissonBC})
directly, we define $\bu = -\bK \nabla p$ and consider the
{\em mixed formulation} of this problem:
\begin{align}
  \bK^{-1} \bu + \nabla p &= 0 \textrm{ in } \Omega \\
  -\nabla \cdot \bu &= -f \textrm{ in } \Omega \\
  p &= g \textrm{ on } \partial \Omega
\end{align}
This problem is discretized using Raviart-Thomas elements, resulting in a
linear system of the form
\begin{equation}
  \left( \begin{array}{cc}
    \bM & \bB^{T} \\ \bB & \bzero
  \end{array} \right)
  \left( \begin{array}{c}
    \bU \\ \bP
  \end{array} \right)
  =
  \left( \begin{array}{c}
    \bbf \\ \bg
  \end{array} \right)
  \label{eq:poissonLinear}
\end{equation}
Block elimination of \eqref{eq:poissonLinear} yields the following block
system:
\begin{align}
  \bS \bP &= \bB \bM^{-1} \bbf - \bg \label{eq:mainPoissonEq} \\
  \bM \bU &= \bbf - \bB^{T} \bP
\end{align}
where $\bS$ is the positive definite Schur complement matrix
$\bS = \bB \bM^{-1} \bB$.  \eqref{eq:mainPoissonEq} can be solved via PCG given
an efficient procedure for forming matrix-vector products with $\bS$.  This
in turn requires efficient application of $\bM^{-1}$.  This will also be
accomplished with the conjugate gradient method.  In the case
$\bK(\bx) = const$, linear systems involving $\bM$ turn out to be quite
easy to solve.  In fact, these systems can be solved quickly using standard
conjugate gradients with no preconditioning.  Nevertheless, this problem
provides a useful benchmark for
our solver.  We also consider versions of the problem in which $\bK(\bx)$
is anisotropic and highly inhomogeneous to demonstrate that our solver
can also handle these cases.

\subsection{Comparisons}
\label{SECcomparisons}
In this section we provide timing and memory usage statistics for some of
the features discussed in \S\ref{SECoffdiagcomp}-\ref{SECoptimizations}. \\

\myparagraph{Diagonal block coordinates:} In \S\ref{SECdiagcoordinates} we
discussed how spatial coordinates for diagonal block indices can be chosen
either geometrically based on information obtained directly from a problem's
discretization, or algebraically via approximations of the low-order
eigenvectors associated with the problem's graph Laplacian.  We compare
both approaches applied to a linear system taken from the nonlinear elasticity
benchmark problem (\S\ref{sec:nonlineardescription}).  We find that both
approaches produce preconditioners which converge in a comparable number of
iterations.  Recall that these spatial coordinates are used to permute
indices within compressed diagonal blocks.  We also consider results when
randomly permuting these indices to demonstrate the need for an effective
permutation.  In this case we observe a significant increase in the number
\begin{wrapfigure}[20]{r}{2.0in}
\centering
\includegraphics[width=2.0in]{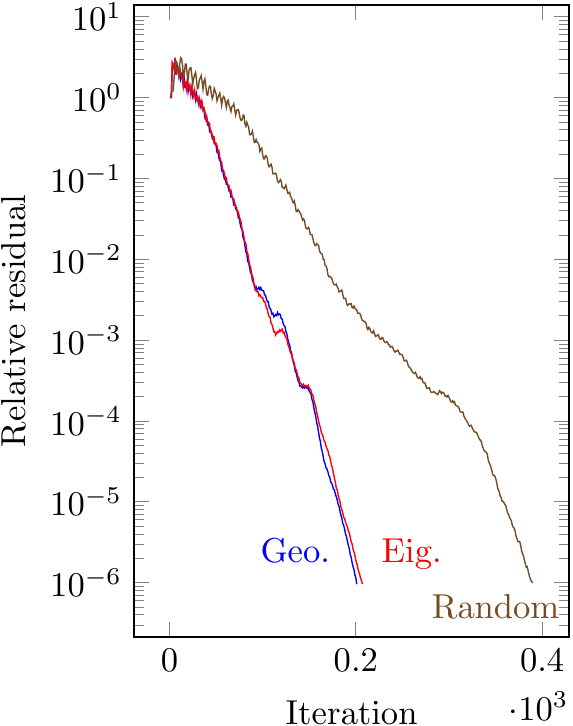}
\figspace{}
\vspace{-4mm}
\caption{PCG convergence using random diagonal block reordering, as well
         as reordering based on geometric index positions (Geo.) and positions
         obtained from a low-accuracy eigensolve (Eig.).
         }
\label{FIGorderingComp}
\end{wrapfigure}
of required PCG iterations.  We also note that the diagonal rank constant
$\tau_{D}$ had to be increased several times in this case to allow for
successful factorization (see \S\ref{SECavoidingIndefinite}).
Convergence plots for these three diagonal block orderings are provided in
Figure~\ref{FIGorderingComp}.
Finally, it
is worth noting that the ordering provided by the original problem discretization
and nested dissection ordering may be sufficient for diagonal block
compression.  For the problem considered here we find that we observe little
difference in the PCG convergence behavior even when no additional permutations
are applied to diagonal block indices.  However, it should be noted that this
is not guaranteed to be the case, and a poor choice of diagonal coordinates
can lead to poor performance (see the random reordering result in
Figure~\ref{FIGorderingComp}). \\

\myparagraph{Interior block performance:} In \S\ref{SECinteriorBlocks} we
discussed how to avoid building and storing certain factor blocks explicitly
to further reduce memory usage.  While this method reduces storage requirements,
it comes at the cost of increased factorization and triangular solve time.  Here
we compare memory usage and solution time results for two examples
computed with and without this optimization.

\section{Conclusions}
\label{SECconclusions}

In this paper, we have described a direct factorization method for the
solution of large sparse linear systems that arise from PDE
discretizations.  Like standard direct solvers, our approach is black
box, and can work with the pre-assembled matrix without prior
information about details of an underlying mesh or a specific
discretization method.  By taking advantage of the low-rank block
structure that arises from the underlying PDE, our method requires
significantly less memory than standard direct methods, but through
careful code organization we retain the high performance of standard
direct solvers through use of level-3 BLAS and LAPACK calls.  We have
demonstrated through examples that our approach retains much of the
robustness of standard direct solvers, and yields a faster time to
solution than the standard multilevel algebraic multigrid
preconditioner ML.

\myparagraph{Limitations and Future Work:}
So far, our work has focused solely on symmetric and positive definite
matrices.  Other authors have showed how to deal with indefinite
problems, with a particular focus on Helmholtz equations, and we
intend to adapt that work along with standard static pivoting
approaches developed in the context of ordinary sparse parallel LU
decompositions.  We also so far only have limited parallelism through
threaded BLAS calls, but intend to extend our code to work in a
distributed memory setting in the future.

\begin{table}[!ht]
\begin{center}
  \resizebox{1.0\hsize}{!}{
  \begin{tabular}{c|p{10cm}}
    \toprule
    Symbol & Meaning \\
    \midrule
    $\bA$
      & System matrix (generally assumed to have been permuted with
        fill-reducing ordering).
      \\
    $\bL$
      & Cholesky factor matrix (of permuted system).
      \\
    $c_{j,s}, c_{j,e}$
      & First and last column in supernode $j$.
      \\
    $\sC_{j}$
      & Supernode $j$'s column set.
      \\
    $\sC_{j}^{O}$
      & Set of columns occurring after supernode $j$.
      \\
    $\sR_{j} = \left\{ r_{j}^{1}, r_{j}^{2}, \ldots \right\}$
      & Set of non-zero rows in supernode $j$'s off-diagonal.
      \\
    $\bL_{j}$
      & Supernode $j$'s block column in $\bL$.
      \\
    $\bLdj, \bLoj$
      & Diagonal and off-diagonal blocks of $\bL_{j}$, respectively.
        $\bLdj \in \bbR^{|\sC_{j}| \times |\sC_{j}|}$ and
        $\bLoj \in \bbR^{|\sR_{j}| \times |\sC_{j}|}$.
      \\
    $\sUdj, \sUdj$
      & Schur complement matrices corresponding to $\bLdj$ and $\bLoj$
      \\
    $\bbD_{j}$
      & Indices of supernode descendants of supernode $j$;
        $\bbD_{j} = \left\{
          1 \le k < j : \sR_{k} \cap \sC_{j} \ne \emptyset \right\}$.
      \\
    $\rDkToj, \rOkToj$
      & Recalling that
        $\sR_{k} = \left\{ r_{k}^{1}, r_{k}^{2}, \ldots \right\}$,
        $\rDkToj = \{ 1 \le p \le |\sR_{k}| : r_{k}^{p} \in \sC_{j} \}$
        and
        $\rOkToj = \{ 1 \le p \le |\sR_{k}| : r_{k}^{p} \in \sR_{j} \}$.
        That is the set of rows in node $k$ needed to construct $\bLdj$ and
        $\bLoj$, respectively (relative to the row space of $\bLok$). \\
    $\srDkToj, \srOkToj$
      & Similar to the definitions above,
        $\srDkToj = \{ r_{k}^{p} \in \sR_{k} : r_{k}^{p} \in \sC_{j} \}$
        and
        $\srOkToj = \{ r_{k}^{p} \in \sR_{k} : r_{k}^{p} \in \sR_{j} \}$.
        That is the set of rows in node $k$ needed to construct $\bLdj$ and
        $\bLoj$, respectively (relative to the full row space of $\bL$). \\
    $\rDkTojs, \srDkTojs$
      & Similar to $\rDkToj$ and $\srDkToj$.
        $\rDkTojs = \{ 1 \le p \le |\sR_{k}| : r_{k}^{p} \in \scDs \}$
        and
        $\srDkTojs = \{ r_{k}^{p} \in \sR_{k} : r_{k}^{p} \in \scDs \}$.
        That is, these are rows from descendant $k$ which are relevant to
        the formation of diagonal block $s$ within $\bLdj$. \\
    $\bVj, \bUj$
      & Low-rank representation for supernode $j$'s off-diagonal block
        ($\bLoj \approx \bVj \bUj^{T}$). \\
    $\dBlock{s}$
      & Block $s$ from the diagonal matrix $\bLdj$. \\
    $\vBlock{s}, \uBlock{s}$
      & Low-rank representation for block $s$ of $\bLdj$. \\
    $\rDs, \cDs$
      & Row and column sets over which $\dBlock{s}$ is defined, relative to
        the dense matrix $\bLdj$. \\
    $\srDs, \scDs$
      & Row and column sets over which $\dBlock{s}$ is defined, relative to
        the full system $\bA$ (or $\bL$).
  \end{tabular}
  }
  \caption {{\bf List of symbols}}
  \label{tab:symbols}
\end{center}
\end{table}

\bibliographystyle{siam}
\bibliography{paper}

\end{document}